\documentclass[a4paper,11pt]{article}
\usepackage{fullpage}
\usepackage{amsmath}
\usepackage{amssymb}
\usepackage{siunitx}

\usepackage[pdftex]{graphicx}
\usepackage[usenames, dvipsnames]{color}

\usepackage{booktabs}
\usepackage{multirow}

\usepackage{subcaption}
\usepackage{authblk}

\usepackage{palatino}
\usepackage{textcomp}

\usepackage{hyperref}

\hypersetup{
bookmarksnumbered,
pdfstartview={FitH},
linkcolor=blue,
citecolor=blue,
colorlinks=true,
pdfpagemode={UseNone},
plainpages=false
}%
\usepackage{soul}

\newcommand{\cV}{\mathcal{V}}

\newcommand{\cE}{\mathcal{E}}
\newcommand{\cK}{\mathcal{K}}

\newcommand{\cU}{\mathcal{U}}

\newcommand{\cP}{\mathcal{P}}

\newcommand{\cvs}{$\text{CVaR}_{0.75}$}
\newcommand{\cvn}{$\text{CVaR}_{0.95}$}

\allowdisplaybreaks

\title{A Comparison of Discrete and Polyhedral Uncertainty Sets for Robust Network Design}

\author[1]{Francis Garuba}
\author[2]{Marc Goerigk\thanks{Corresponding author. Email: \texttt{marc.goerigk@uni-siegen.de}}}
\author[3]{Peter Jacko}

\affil[1]{{\small Mathematical Sciences, University of Southampton, Highfield, Southampton, SO17 1BJ, United Kingdom}}
\affil[2]{{\small Department of Management Science, Lancaster University, Lancaster LA1 4YX, United Kingdom}}
\affil[3]{{\small Network and Data Science Management, University of Siegen, Unteres Schlo{\ss} 3, 57072 Siegen, Germany}}

\date{}

\begin{document}

\maketitle

\begin{abstract}
We consider a network design and expansion problem, where we need to make a capacity investment now, such that uncertain future demand can be satisfied as closely as possible. To use a robust optimization approach, we need to construct an uncertainty set that contains all scenarios that we believe to be possible. 
In this paper we discuss how to construct two common types of uncertainty sets, which are discrete and polyhedral uncertainty, using real-world data. We employ clustering to generate a discrete uncertainty set, and place hyperplanes sequentially to generate a polyhedral uncertainty set. 
We then compare the performance of the resulting robust solutions for these two types of models on real-world data. Our results indicate that polyhedral models, while being popular in the recent network design literature, are less effective than discrete models both in terms of computational burden and solution quality with respect to the performance measure considered.
\end{abstract}


\section{Introduction}

Operations research approaches have found wide application in the planning, design and operations management of transportation, power and energy distribution, supply chain logistics and telecommunications networks. In particular, many types of optimization models have been developed over the last decades for network design and expansion problems, see, e.g., \cite{Bertsekas1998,Magnanti1984,Minoux1989}. Optimization of investments has thus attained a key strategic role in this industry. Moreover, these decisions need to be made well ahead of time based on a forecast of future traffic demand.
Unfortunately, traffic demand has proven to be difficult to predict accurately. In order to factor in this uncertainty and design a network that is immune to traffic variability, robust optimization approaches have been proposed. For this purpose, a number of uncertainty models have already been developed and
investigated (see \cite{Ben-Tal2009,Bertsimas2011,goerigk2016algorithm}). The drawback of classic approaches, however, is that the uncertainty set is assumed to be given, i.e., the decision maker can advise on how the uncertainty is shaped. Moreover, an inappropriate choice of uncertainty set may result
in models that are too conservative or in some cases computationally intractable. As the decision maker cannot be expected to make this choice in practice, data-driven and learning approaches have been proposed as well (see \cite{Bertsimas2017,Chassein2018}).

Robust optimization in general has found increasing use and application in the network design area. In \cite{Atamt_rk_2007} a two-stage robust
network flow problem under demand uncertainty was considered following the work of \cite{Ben_Tal_2004}, while affine routing in the their
robust network capacity planning model was used in \cite{Ouorou2007}. \cite{Ord_ez_2007} looked at network capacity expansion under both demand and cost uncertainty, while \cite{Koster_2013}
considered a robust network design problem with static routing in the setting of \cite{Bertsimas_2004}. \cite{Poss_2012} considered robust network design
with polyhedral uncertainty and \cite{Babonneau2013} robust capacity assignment for networks with uncertain demand. \cite{Pessoa_2015} used a cutting plane
algorithm while taking into consideration the uncertainty in unmet demand outsourced cost.
Regarding uncertainty sets, polyhedral models are most frequently used in radio network design, along with hose models from the works of
\cite{Duffield1999,Fingerhut_1997}, budget uncertainty by \cite{Atamt_rk_2007}, cardinal constrained uncertainty by \cite{Bertsimas_2004}, and interval
uncertainty, among others. Little research compares these models of uncertainty. \cite{Atamt_rk_2007} compared their single-stage robust model
using budget uncertainty with a scenario-based two-stage stochastic approach. \cite{Chassein2018} constructed different uncertainty sets from real world
data and compared performance within and outside sample for shortest path problems.

The basic network design problem that we consider in this paper is as follows. Given an undirected graph $G=(\cV,\cE)$ and currently installed
capacity $u_e$ for each edge $e\in \cE$, we would like to determine an amount of capacity $x_e$ to be installed additionally. For each edge $e$, we are
given an investment cost $c_e$ per unit of additionally installed capacity. As the graph is undirected, the direction of flow is not relevant for our
model, and we define $\cK = \{ \{i,j\} : i,j \in \cV, i < j\}$ as the set of commodities, where each commodity $ k $ is identified by an unordered
pair of nodes $\{i,j\}$ between which a given demand needs to be satisfied. Let $d_k$ be the demand corresponding to commodity $k\in\cK$, and let $\cP_k$
be the set of simple paths in $G$ connecting the nodes of the commodity. The aim is to find capacities $\pmb{x}$ such that all demands are fulfilled and
the capacity expansion costs are as small as possible. Formally, the baseline model can thus be written as follows.
\begin{subequations}
\begin{align}
\min\ & \sum_{e\in \cE} c_e x_e \\
\text{s.t. } & \sum_{p\in \cP_k} f_{kp} \ge d_k & \forall k\in\cK \label{nom1}\\
& \sum_{k\in\cK} \sum_{p\in \cP_k:e\in p} f_{kp} \le u_e + x_e & \forall e\in\cE \label{nom2}\\
& x_e \ge 0 & \forall e\in\cE \\
& f_{kp} \ge 0 & \forall k\in\cK, p\in \cP_k.
\end{align}
\end{subequations}
The variables $f_{kp}$ model the amount of flow along path $p$ for commodity $k$. Here, Constraints~\eqref{nom1} ensure that a sufficient amount of flow is
transported along all paths connecting source and sink of commodity $k\in\cK$, while Constraints~\eqref{nom2} model that each edge needs to provide
sufficient capacity. Instead of using a path-based formulation, it is also possible to use a model with flow variables for every edge in the network
(see, e.g., \cite{Gendron1999}). In this paper, we focus on the path-based formulation, as it performed better in our computational experiments.

In practice, the demand $\pmb{d}$ changes over time and is not known precisely. Thus, a two-stage model is required, where we decide now where to build how much capacity (the strategic decision $\pmb{x}$), and we can decide where to route the flow once the demand is known (the operational decision $\pmb{f}$). Let us assume that a set $\cU$ can be constructed that contains all demand scenarios $\pmb{d}$ that we would like to take into account for our planning. The two-stage robust network design problem is then to solve
\begin{subequations}
\label{two-stage-model}
\begin{align}
\min\ & \sum_{e\in \cE} c_e x_e \label{mod1}\\
\text{s.t. } & \sum_{p\in \cP_k} f_{kp}(\pmb{d}) \ge d_k & \forall k\in\cK, \pmb{d}\in\cU \label{mod2}\\
& \sum_{k\in\cK} \sum_{p\in \cP_k:e\in p} f_{kp}(\pmb{d}) \le u_e + x_e & \forall e\in\cE,\pmb{d}\in\cU \label{mod3}\\
& x_e \ge 0 & \forall e\in\cE \label{mod4}\\
& f_{kp}(\pmb{d}) \ge 0 & \forall k\in\cK, p\in \cP_k, \pmb{d}\in\cU. \label{mod5}
\end{align}
\end{subequations}
In this setting, $f_{kp}$ has become a function that depends on the scenario $\pmb{d}$. Note that in Constraint~\eqref{mod2}, $ d_k $ is a component of $\pmb{d}$, thus is also scenario-dependent.

In model~\eqref{two-stage-model}, there is one set of variables $\pmb{f}$ for each possible scenario $\pmb{d}\in\cU$. If $\cU$ is a discrete set, then the corresponding model remains a linear program with $|\cE| + |\cU|\cdot \sum_{k\in\cK} |\cP_k|$ many variables. Hence, smaller sets $\cU$ result in models that are easier to solve. For continuous sets $\cU$ such as polyhedra, model~\eqref{two-stage-model} has infinitely many constraints and variables and cannot be solved directly. The most common method of approaching such problems is to use the affinely adjustable counterpart \cite{Ben_Tal_2004}. Here, the crucial factor determining the resulting problem size is the number of constraints and variables that define the uncertainty polyhedron. Current approaches such as the kernel learning method \cite{ning2018data,shang2017data} or the 
Bayesian nonparametric set construction approach \cite{campbell2015bayesian,ning2017data} create (unions of) polyhedra with tend to have a high number of constraints and variables even for low-dimensional problems. Hence, for the network design problem considered here, we require a new approach to generate polyhedra that give a reasonable description of the uncertainty with a small number of constraints and no auxiliary variables.

In this paper we consider the question whether discrete or polyhedral uncertainty sets are more practicable for robust network design. We use a clustering approach (using the well-known $K$-means clustering data mining method) to generate discrete uncertainty sets and calculate the cluster centroids for real-world data taken from SNDlib (see \cite{Orlowski2010}). We compare this solution to the solution obtained when modelling uncertainty using a polyhedral set, where constraints on the demand are given by hyperplanes that are generated dynamically. To the best of our knowledge, this is for the first time that such a comparison is done for network design problems. For the real-world dataset we consider, we find in our numerical experiments that solutions based on discrete uncertainty outperform solutions based on polyhedral uncertainty when using high risk-adverse performance metrics such as maximum or {\cvn} of unsatisfied demand. At the same time, solutions based on discrete uncertainty found by clustering can be computed two orders of magnitude faster than those based on polyhedral uncertainty.

The rest of this paper is organized as follows. As the problem data is the center point of our research, we first discuss this in Section~\ref{sec:data}.
In particular, we describe how to construct uncertainty sets $\cU$ from the data. We then introduce models for robust network design for both discrete and
polyhedral uncertainty sets in Section~\ref{sec:models}. Experimental results are discussed in Section~\ref{sec:experiments}. Finally,
Section~\ref{sec:conclusions} concludes our work and points out future research directions.

\section{Problem Data and Uncertainty Set Construction}
\label{sec:data}

\begin{figure}[htbp]
\begin{center}
\includegraphics[width=.55\textwidth]{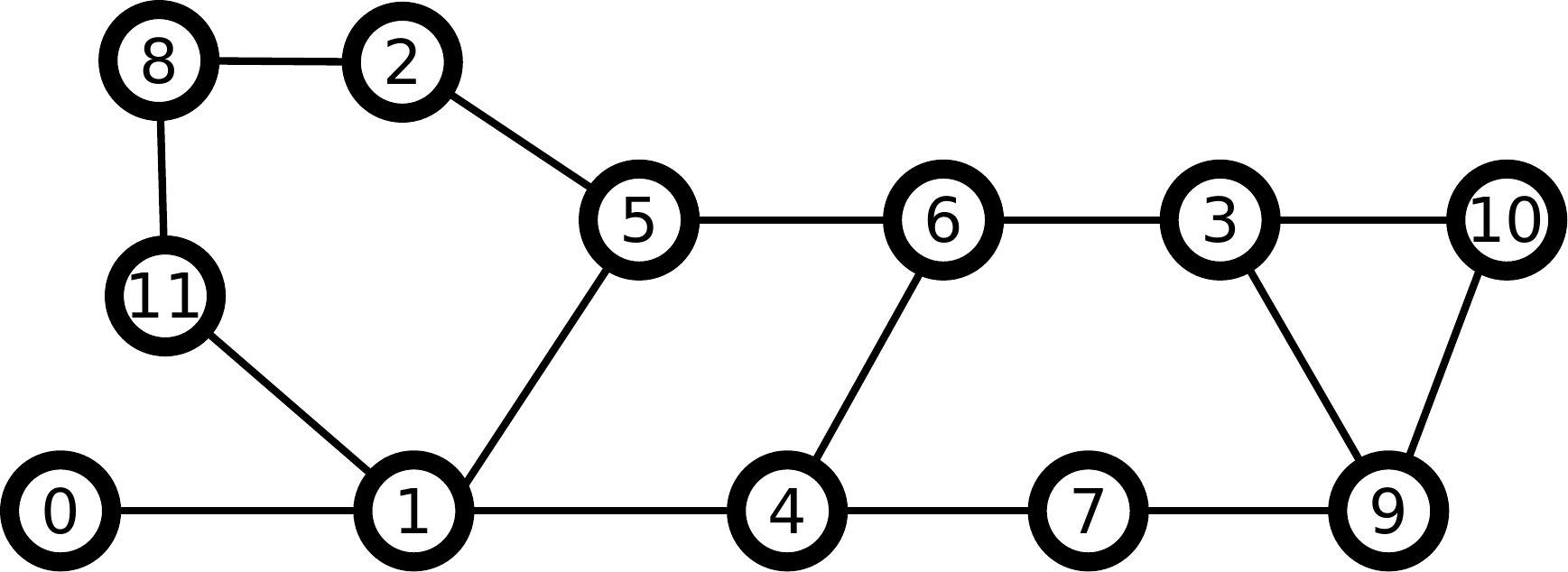}
\caption{Abilene network topology.}\label{fig:abilene}
\end{center}
\end{figure}

In the following discussion, we focus on the Abilene network based on data from the SNDlib (see \cite{Orlowski2010}). It consists of 12 nodes connected by 15 edges, see
Figure~\ref{fig:abilene}, which spread over the US. With 12 nodes, there exist $12\cdot11/2 = 66 =: \kappa$ different possible commodities.

Data was collected by Yin Zhang\footnote{\url{http://www.cs.utexas.edu/~yzhang/}} in 5 minute intervals between 01.03.2004 and 10.09.2004 with some breaks
in between. Table~\ref{tab:scenarios} shows the number of measurements that are available for each month. Note that one day can give 288 measurements in 5
minute intervals. Based on this number, we also show the maximum number of possible measurements that can be achieved each month, but note that not all
data are available.

\begin{table}[htb]
\begin{center}
\begin{tabular}{r|rrrrrrr}
Month & 03 & 04 & 05 & 06 & 07 & 08 & 09 \\
\hline
\# Measurements available & 4,032 & 6,048 & 8,928 & 8,640 & 8,928 & 8,640 & 2,880 \\
\# Measurements possible &  8,928 & 8,640 & 8,928 & 8,640 & 8,928 & 8,928 & 8,640
\end{tabular}
\caption{Numbers of available measurements for each month.}\label{tab:scenarios}
\end{center}
\end{table}

\begin{figure}[htbp]
\begin{center}
\begin{subfigure}{0.49\textwidth}
\centering
\includegraphics[width=\linewidth]{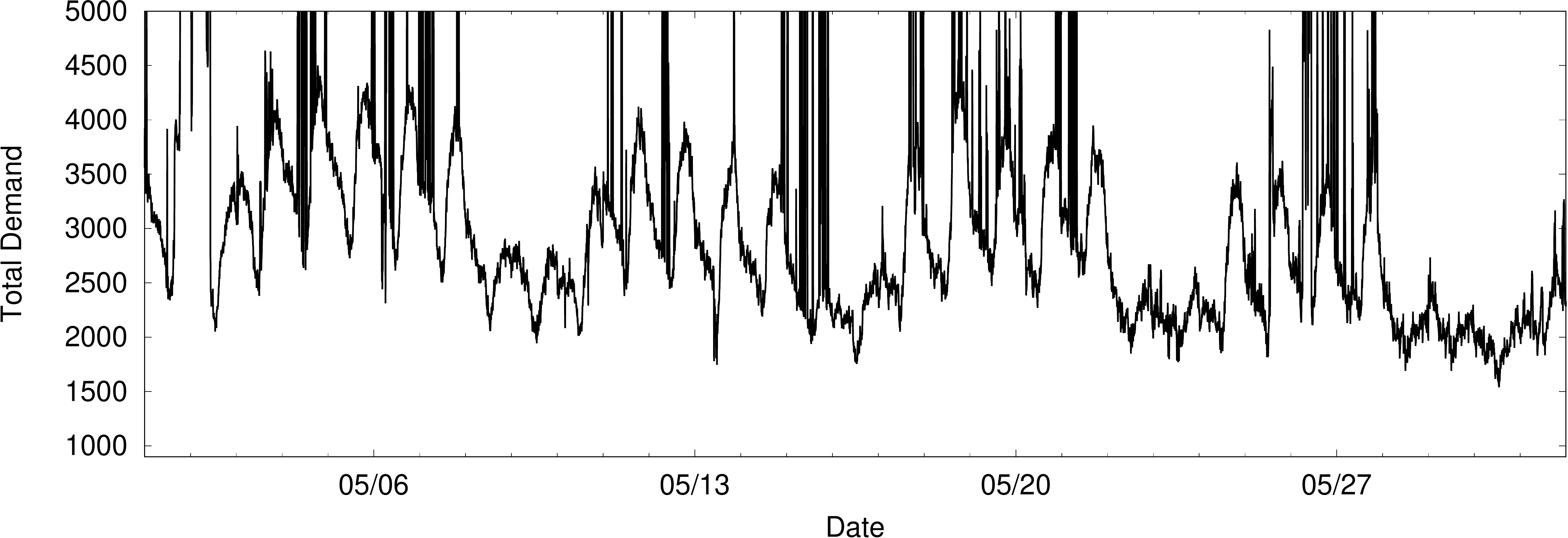}
\caption{Demand profile for month 05.}\label{fig:data05}
\end{subfigure}\hfill
\begin{subfigure}{.49\textwidth}
\centering
\includegraphics[width=\linewidth]{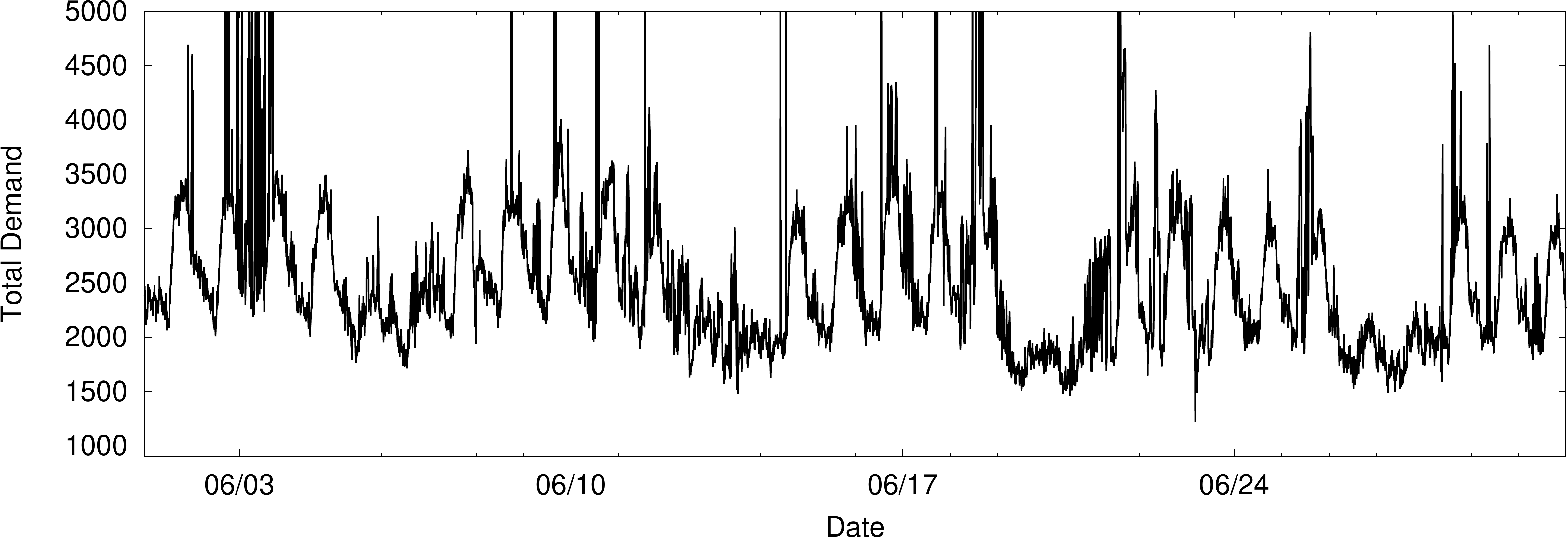}
\caption{Demand profile for month 06.}\label{fig:data06}
\end{subfigure}
\begin{subfigure}{.49\textwidth}
\centering
\includegraphics[width=\linewidth]{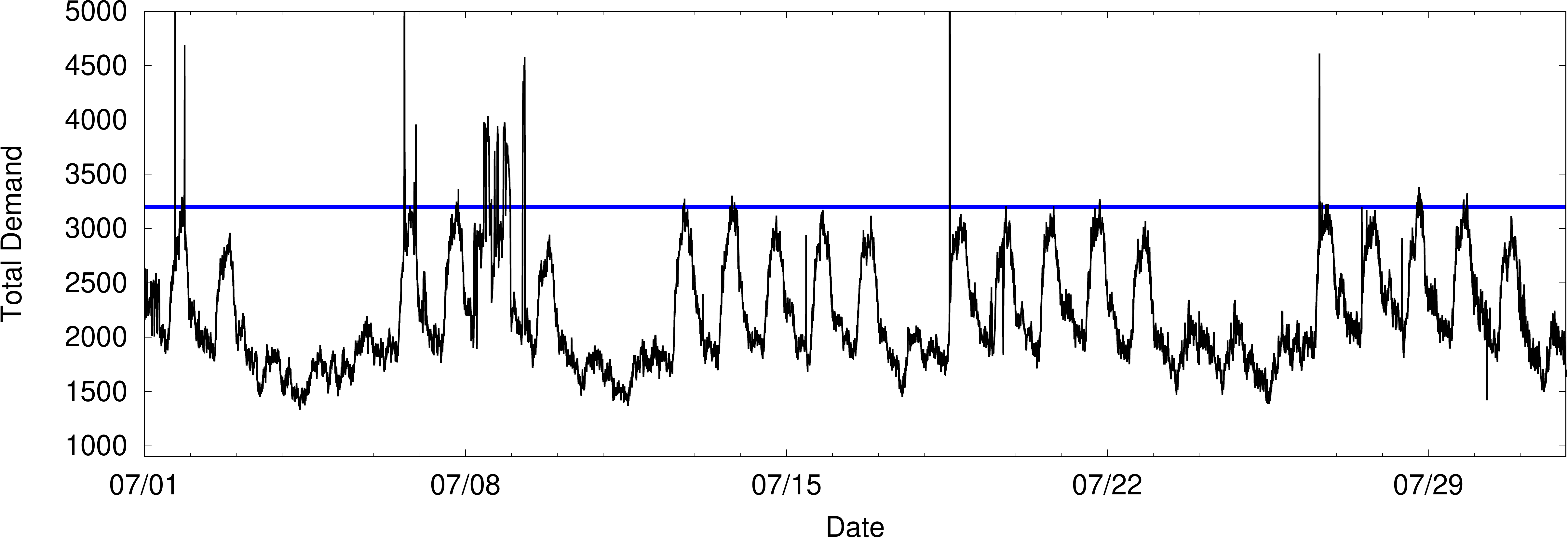}
\caption{Demand profile for month 07 (training set).}\label{fig:data07}
\end{subfigure}\hfill
\begin{subfigure}{.49\textwidth}
\centering
\includegraphics[width=\linewidth]{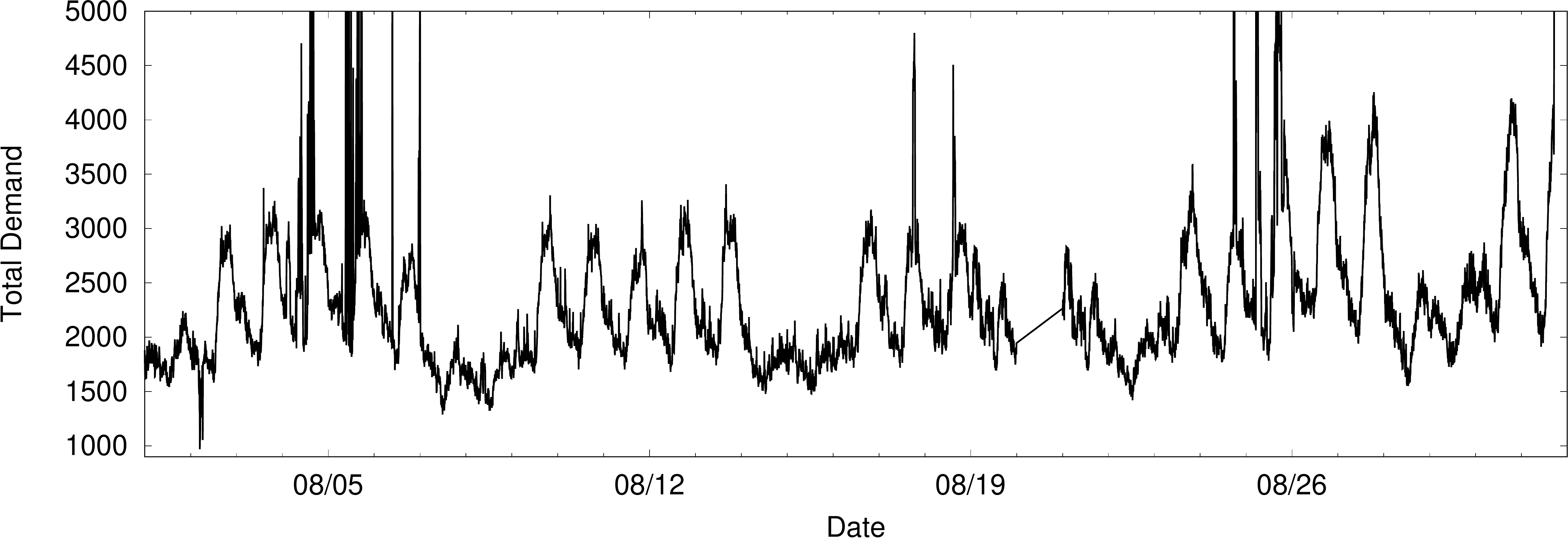}
\caption{Demand profile for month 08.}\label{fig:data08}
\end{subfigure}
\caption{Total demand in the Abilene network.}\label{fig:alldata}
\end{center}
\end{figure}

We focus on months 05--08, where data are most complete, and use an arbitrary month for training (in this case, month 07), and the remaining data for testing. The total demand per scenario, i.e., the sum of demand over all commodities, for months 05--08 is presented in Figure~\ref{fig:alldata}. The horizontal line indicates a cut-off value at the 98\%-quantile, where we ignore data with total demand above the line for training.

Let $T$ be the number of measurements available for training, each being a $ \kappa $-dimensional vector of reals. 
We now discuss how to generate discrete and polyhedral uncertainty sets based on the training data points. Let $\mathcal{D} =
\{\pmb{d}^1,\ldots,\pmb{d}^T\}$ denote this training set. For a discrete uncertainty set $\cU^d$, where each scenario is explicitly listed, a natural
approach is setting $\cU^d = \mathcal{D}$. But it has been shown (see \cite{Chassein2018}) that this can lead to an overfitting effect, such that the
resulting robust solutions do not perform well on out-of-sample data points. Furthermore, it is desirable to control the degree of conservatism. We
therefore propose a clustering approach to generate discrete uncertainty sets. We aggregate similar scenarios together, with the intention to reduce the
problem complexity on the one hand, and to become less dependent on data noise on the other hand, as outliers are reduced by merging them with other scenarios.
Scenario aggregation based on $K$-means clustering has
been applied as an approximation method also to robust min-knapsack problems, see \cite{chassein2018approximating}. Let $\cU^d_K$ denote a discrete
uncertainty set derived from a $K$-means clustering of the set $\mathcal{D}$. Then these sets contain the two special cases $\cU^d_T = \mathcal{D}$, i.e., the
original set of training points, and $\cU^d_1$, which consists of only the average case scenario.

\begin{figure}[htbp]
\begin{center}
\begin{subfigure}{.4\textwidth}
\includegraphics[width=\linewidth]{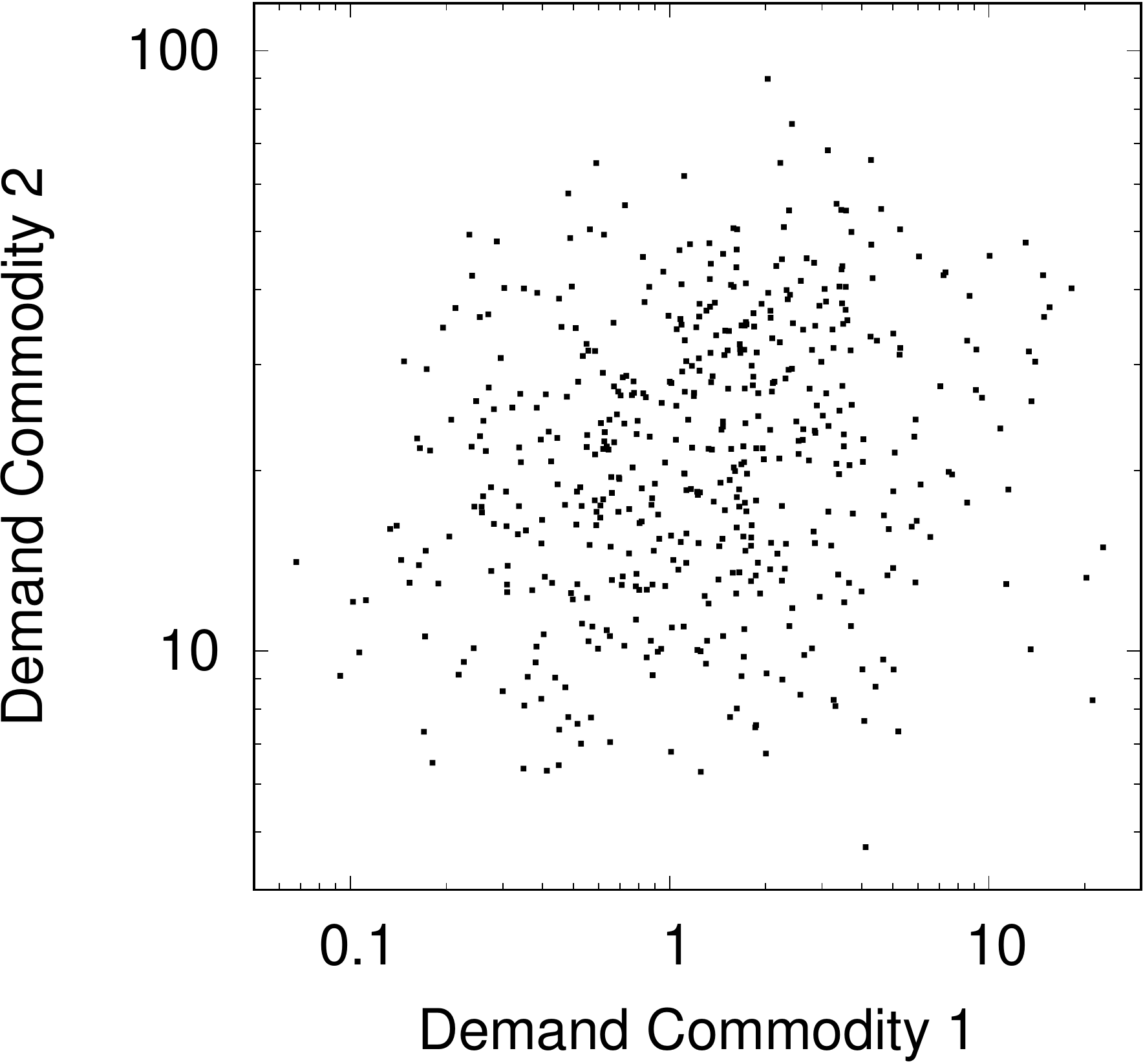}
\caption{Subset of training data $\mathcal{D}$.}\label{fig:dataa}
\end{subfigure}
\hspace{1cm}
\begin{subfigure}{.4\textwidth}
\includegraphics[width=\linewidth]{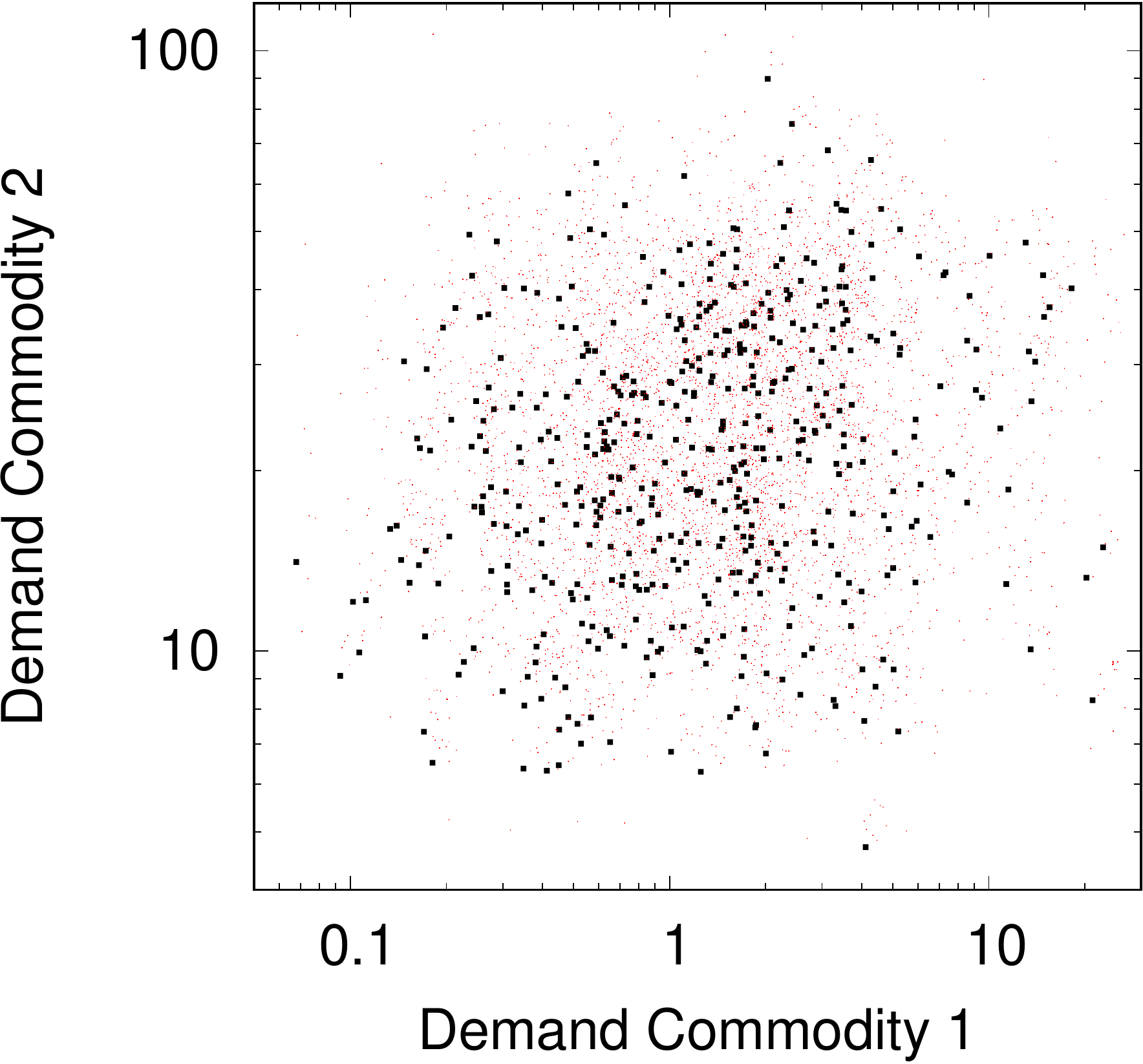}
\caption{Training data with added noise data in red.}\label{fig:datab}
\end{subfigure}

\vspace*{5mm}

\begin{subfigure}{.4\textwidth}
\includegraphics[width=\linewidth]{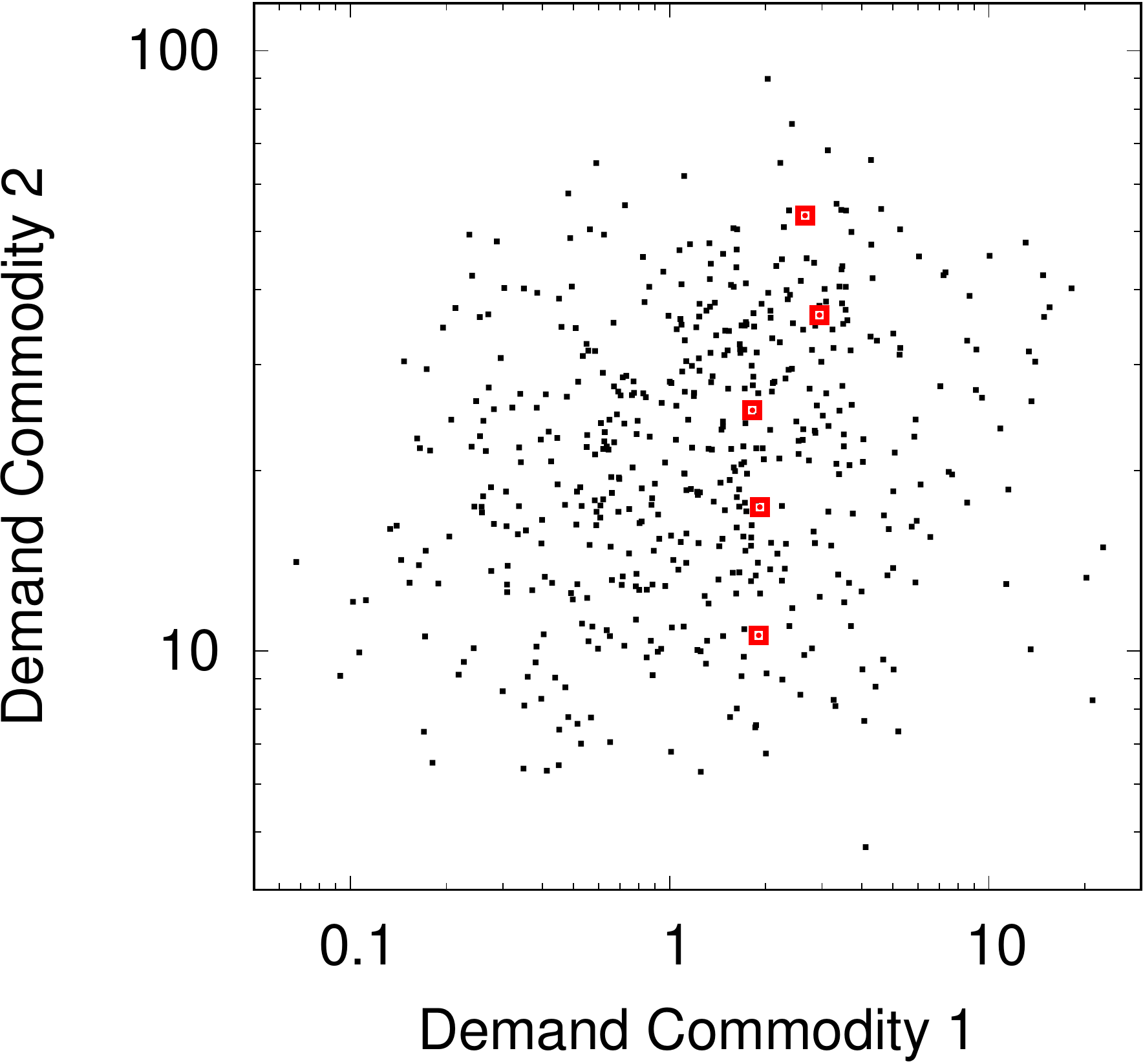}
\caption{A $K$-means clustering with $K=5$ in red.}\label{fig:datac}
\end{subfigure}
\hspace{1cm}
\begin{subfigure}{.4\textwidth}
\includegraphics[width=\linewidth]{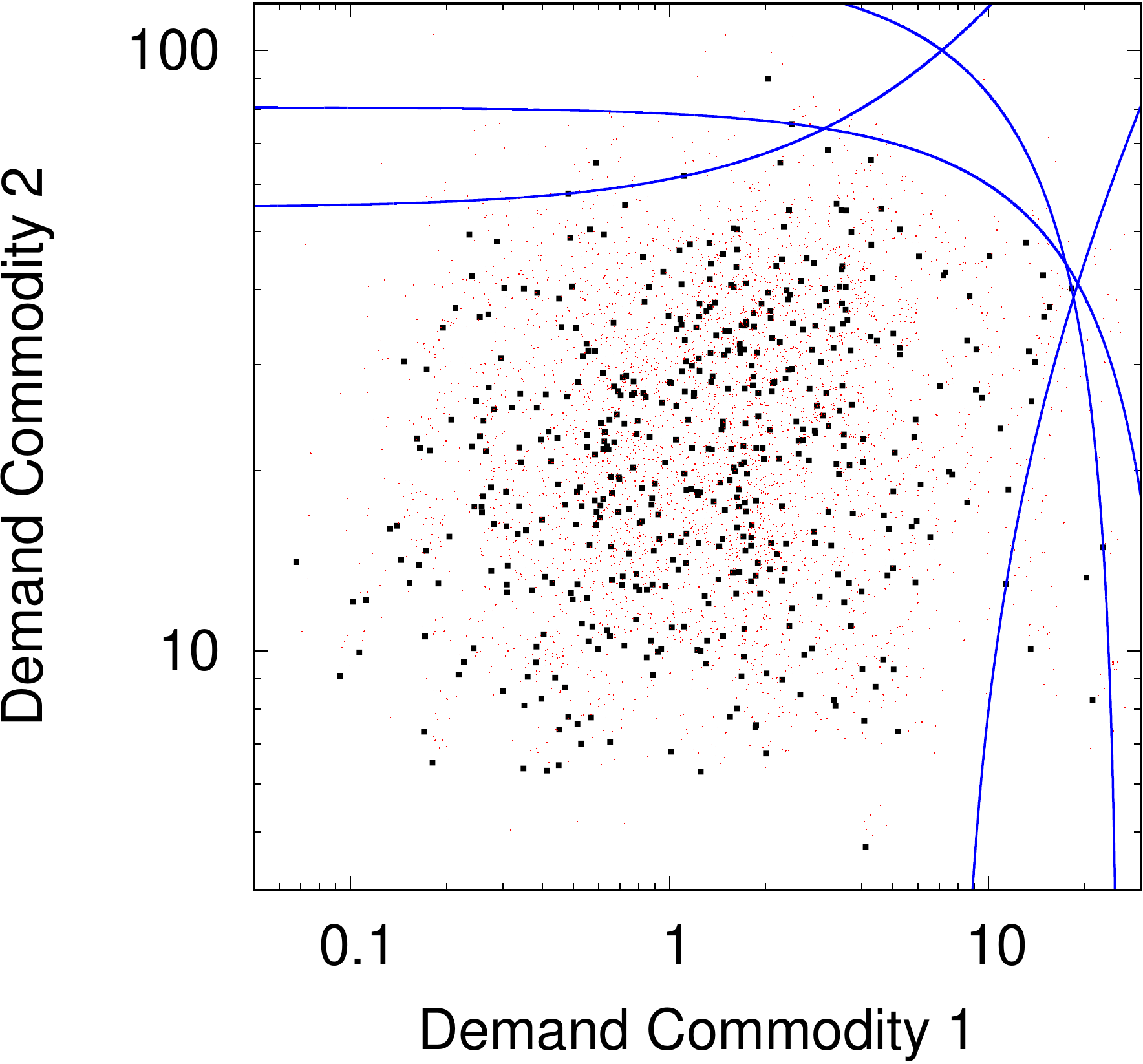}
\caption{Polyhedron with $M=4$ hyperplanes.}\label{fig:datad}
\end{subfigure}

\vspace*{5mm}

\begin{subfigure}{.4\textwidth}
\includegraphics[width=\linewidth]{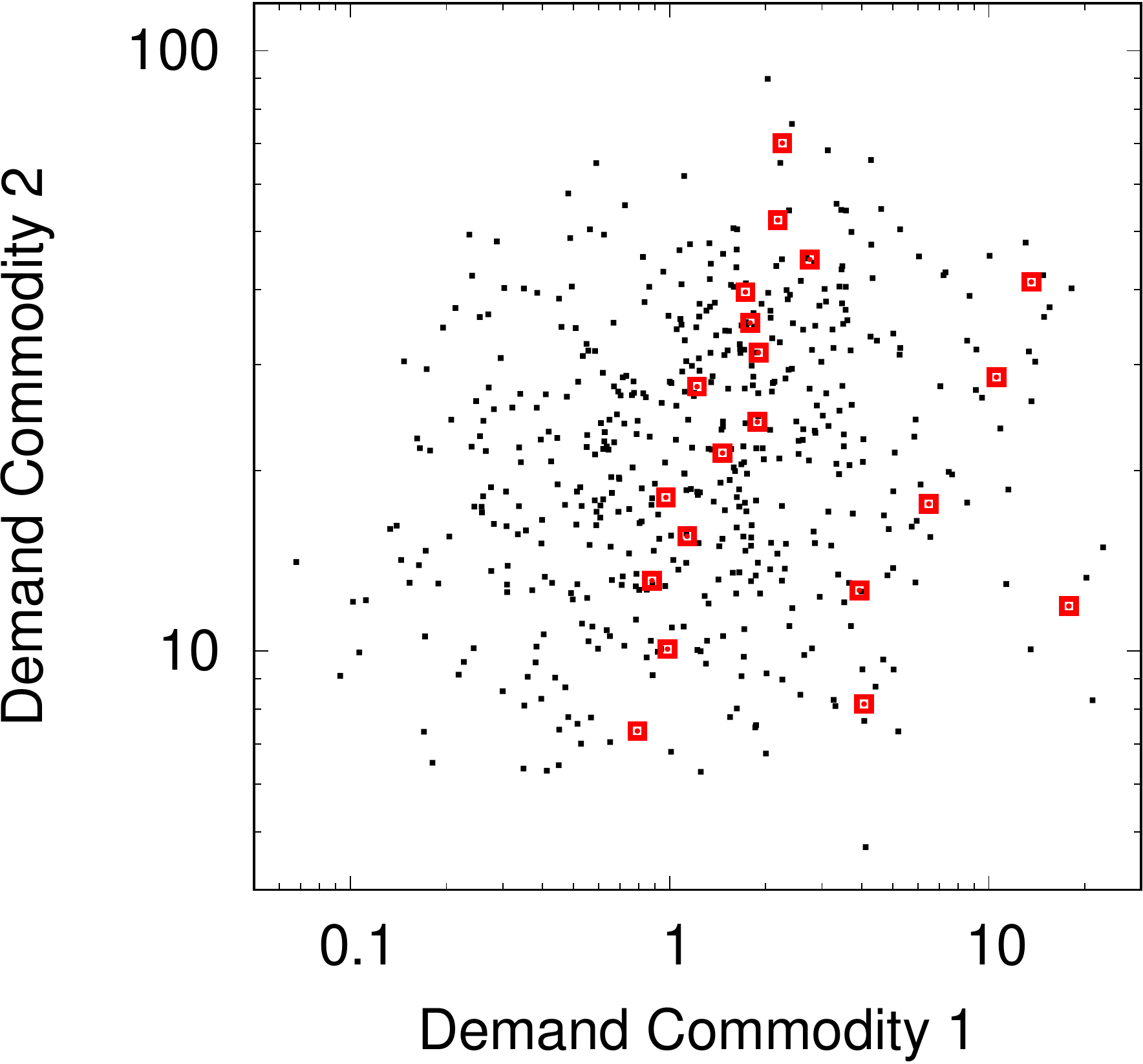}
\caption{A $K$-means clustering with $K=20$ in red.}\label{fig:datae}
\end{subfigure}
\hspace{1cm}
\begin{subfigure}{.4\textwidth}
\includegraphics[width=\linewidth]{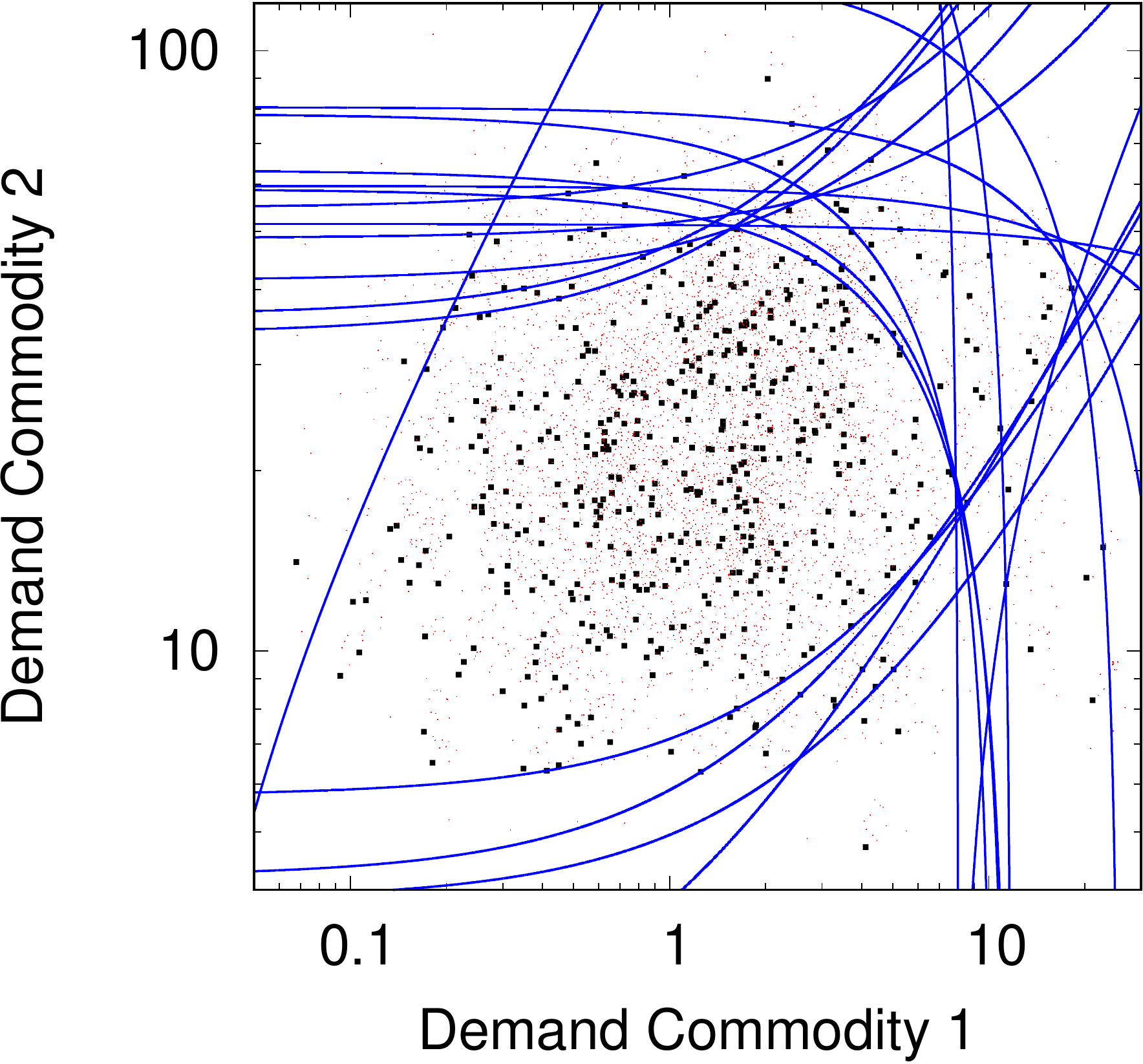}
\caption{Polyhedron with $M=20$ hyperplanes.}\label{fig:dataf}
\end{subfigure}
\caption{Illustration of methods to generate discrete and polyhedral uncertainty on a subset of training data restricted to two dimensions (note the logarithmic scale of the axes).}\label{fig3}
\end{center}
\end{figure}

We show a simple example in Figure~\ref{fig3}. In Figure~\ref{fig:dataa}, we plot a subset of the training data, restricted to two arbitrarily chosen
dimensions (recall that every commodity corresponds to a dimension of the demand vector). In Figure~\ref{fig:datac}, we show a discrete uncertainty
$\cU^d_5$ set based on a $K$-means clustering with $K=5$ centers, which captures the training data only in a rough manner. With $K=20$ (see
Figure~\ref{fig:datae}), most features of the data have been captured.

We now consider the case of polyhedral uncertainty,
\[ \cU^p_M = \big\{ \pmb{d}\in\mathbb{R}^{\kappa}_+ : V\pmb{d} \le \pmb{b}, d_k\in[\underline{d}_k,\overline{d}_k] \big\} \]
where $V=(v_{ik})$ is a matrix in $\mathbb{R}^{M\times \kappa}$ and $\pmb{b}$ is a vector in $\mathbb{R}^{\kappa}$ (i.e., there are $M$ linear constraints
on the demand vector). As the number of constraints $M$ has a significant impact on the solution time of the resulting robust model due to its two-stage nature, we would like to find only few constraints which describe the training data $\mathcal{D}$ well. 
This excludes the application of existing methods that generate polyhedra with high number of constraints or additional dimensions \cite{shang2017data}. The approach we follow here is to place hyperplanes one by one as efficiently as possible. To this end, we generate a set of noise data points, which we would like to distinguish from the original training data by placing hyperplanes that
put as many original points on one side, and as many noise points on the other side as possible. This trade-off is adjusted dynamically: for the first
hyperplane, there is a high penalty for original points that are classified as noise. This way, we find an outer description of the data, which results in
large and conservative uncertainty sets. This penalty is reduced over time, so that later hyperplanes become less conservative and cut away more outliers in the
training data. Noise points are generated by randomly increasing values of single training data points, and randomly using values from other data points in
single dimensions with low probability.

In Figure~\ref{fig3}, we use the same data as for the clustering example to illustrate this process. The random noise is shown as red points in Figure~\ref{fig:datab}. The first four hyperplanes we generate are given in Figure~\ref{fig:datad}. It can be seen that they form an outer approximation of the data, removing only few outliers in the process. With an increasing number of hyperplanes $M$, the polyhedron $\cU^p_M$ becomes smaller and less conservative (see Figure~\ref{fig:dataf}).

\section{Robust Models}
\label{sec:models}

We now discuss how to reformulate the general model~\eqref{two-stage-model} for specific uncertainty sets $\cU^d_K$ and $\cU^p_M$.

\subsection{Discrete Uncertainty}

Let $\cU^d_K = \{\pmb{d}^1,\ldots,\pmb{d}^K\}$ be given. As this set is discrete, we can simply write $f_{kp}(\pmb{d}^i) = f^i_{kp}$ for all $\pmb{d}^i\in\cU^d_K$. We write $[K]:=\{1,\ldots,K\}$ in the following. The resulting compact optimization model is then given as follows:
\begin{subequations}
\label{disc-model}
\begin{align}
\min\ & \sum_{e\in \cE} c_e x_e \label{disc-start}\\
\text{s.t. } & \sum_{p\in \cP_k} f^i_{kp} \ge d^i_k & \forall k\in\cK, i\in[K] \\
& \sum_{k\in\cK} \sum_{p\in \cP_k:e\in p} f^i_{kp} \le u_e + x_e & \forall e\in\cE, i\in[K] \\
& x_e \ge 0 & \forall e\in\cE \\
& f^i_{kp} \ge 0 & \forall k\in\cK, p\in \cP_k, i\in[K]. \label{disc-end}
\end{align}
\end{subequations}
Model~\eqref{disc-model} is a linear program and can be solved directly using standard methods.

\subsection{Polyhedral Uncertainty}

Rewriting $f_{kp}(\pmb{d})$ in a compact form is less straightforward for continuous uncertainty sets than in the previous discrete case. We apply the well-known affine decision rules (also known as affine adjustable robust counterpart) approach, see \cite{Ben_Tal_2004}. To this end, we restrict $f_{kp}(\pmb{d})$ to be an affine linear function in $\pmb{d}$ by writing $f_{kp}(\pmb{d}) = \phi_{kp} + \sum_{\ell\in\cK} d_\ell \Phi_{kp\ell}$. 
Here, $\phi_{kp} \gtrless 0$ and $\Phi_{kpl} \gtrless 0$ are new decision variables. By using affine decision rules, we restrict the set of feasible solutions, and thus form a conservative approximation to the original problem. The advantage is that second-stage decision variables $f_{kp}(\pmb{d})$ are expressed through first-stage decision variables, which means that the min-max-min problem is reduced to a min-max problem.

By substituting the $f_{kp}(\pmb{d})$ variables in \eqref{two-stage-model} and rearranging terms,
this problem becomes:
\begin{subequations}
\begin{align}
\min\ & \sum_{e\in \cE} c_e x_e \\
\text{s.t. } & \sum_{p\in \cP_k} \phi_{kp} \ge \max_{\pmb{d}\in\cU^p_M} \sum_{\ell\in\cK} \left( 1_{\ell=k} - \sum_{p\in \cP_k} \Phi_{kp\ell}\right) d_\ell & \forall k\in\cK \label{con17}\\
& \sum_{k\in\cK} \sum_{p\in \cP_k : e\in p} \phi_{kp} + \max_{\pmb{d}\in\cU^p_M} \sum_{\ell\in\cK} \left( \sum_{k\in\cK} \sum_{p\in \cP_k : e\in p}  \Phi_{kp\ell}\right) d^e_\ell  \le u_e + x_e \hspace{-1cm} & \forall e \in\cE\\
& \phi_{kp} + \min_{\pmb{d}\in\cU^p_M} \sum_{\ell\in\cK} \Phi_{kp\ell} d_{\ell} \ge 0 & \forall k\in\cK,p\in \cP_k \\
& x_e \ge 0 & \forall e\in\cE \\
& \phi_{kp} \gtrless 0 & \forall k\in\cK, p\in \cP_k \\
& \Phi_{kpl} \gtrless 0 & \forall k\in\cK, p\in \cP_k,\ell\in\cK.
\end{align}
\end{subequations}
Here, $1_{\ell=k}$ denotes the indicator function, which is one if $\ell=k$, and 0 otherwise.
The inner maximization and minimization problems can then be reformulated using linear programming duality. As an example, consider Constraint~\eqref{con17} for a fixed $k\in\cK$. The value of the right-hand side is
\begin{equation}
\label{model-prim}
\max\left\{ \sum_{\ell\in\cK} \left( 1_{\ell=k} - \sum_{p\in \cP_k} \Phi_{kp\ell}\right) d_\ell :   V\pmb{d} \le \pmb{b},\ d_\ell \in [\underline{d}_\ell, \overline{d}_{\ell}] \  \forall \ell\in\cK\right\}
\end{equation}
As $\cU^p_M$ is polyhedral, this is a linear program, the dual of which is
\begin{subequations}
\label{model-dual}
\begin{align}
\min\ & \sum_{i\in[M]} b_i \alpha_{i} + \sum_{\ell\in\cK} (\overline{d}_\ell\overline{\beta}_{\ell} - \underline{d}_\ell \underline{\beta}_{\ell}) \label{dual1}\\
\text{s.t. } & \sum_{i\in[M]} v_{i\ell} \alpha_{i} + \overline{\beta}_{\ell} - \underline{\beta}_{\ell} \ge 1_{\ell=k} - \sum_{p\in \cP_k} \Phi_{kp\ell} & \forall \ell\in\cK \label{dual2}\\
& \pmb{\alpha},\overline{\pmb{\beta}}, \underline{\pmb{\beta}} \ge 0. \label{dual3}
\end{align}
\end{subequations}
By weak duality, any feasible solution to~\eqref{model-dual} gives an upper bound on the value of~\eqref{model-prim}. By strong duality, the best such bound is equal to the optimal primal objective value.
Thus we can substitute the formulation~\eqref{model-dual} into Constraint~\eqref{con17} to reach an equivalent linear reformulation. Repeating this process for all constraints, the robust network extension problem with polyhedral uncertainty can be rewritten in the following way:
\begin{align*}
\min\ &\sum_{e\in\cE} c_e x_e \\
\text{s.t. } & \sum_{p\in \cP_k} \phi_{kp} \ge \sum_{i\in[M]} b_i \alpha_{ki} + \sum_{\ell\in\cK} (\overline{d}_\ell\overline{\beta}_{k\ell} - \underline{d}_\ell \underline{\beta}_{k\ell}) & \forall k\in\cK \\
 & \sum_{i\in[M]} v_{i\ell} \alpha_{ki} + \overline{\beta}_{k\ell} - \underline{\beta}_{k\ell} \ge 1_{\ell=k} - \sum_{p\in \cP_k} \Phi_{kp\ell} & \forall k,\ell\in\cK \\
&\sum_{k\in\cK} \sum_{p\in \cP_k : e\in p} \phi_{kp} + \sum_{i\in[M]} b_i \pi_{ei} + \sum_{\ell\in\cK} (\overline{d}_\ell\overline{\rho}_{e\ell} - \underline{d}_\ell\underline{\rho}_{e\ell} )
\le u_e + x_e & \forall e \in\cE \\
& \sum_{i\in[M]} v_{i\ell} \pi_{ei} + \overline{\rho}_{e\ell} - \underline{\rho}_{e\ell} \ge \sum_{k\in\cK} \sum_{p\in \cP_k : e\in p}  \Phi_{kp\ell} & \forall e\in\cE, \ell\in\cK\\
& \phi^k_p \ge \sum_{i\in[M]} b_i\xi_{kpi} + \sum_{\ell\in\cK} (\overline{d}_\ell \overline{\zeta}_{kp\ell} - \underline{d}_\ell \underline{\zeta}_{kp\ell}) & \forall k\in\cK,p\in \cP_k \\
& \sum_{i\in[M]} v_{i\ell} \xi_{kpi} + \overline{\zeta}_{kp\ell} - \underline{\zeta}_{kp\ell} \ge - \Phi_{kp\ell} & \forall k,\ell\in\cK,p\in \cP_k \\
& x_e \ge 0 & \forall e\in\cE \\
& \phi_{kp} \gtrless 0 & \forall k\in\cK, p\in \cP_k \\
& \Phi_{kp\ell} \gtrless 0 & \forall k\in\cK, p\in \cP_k,\ell\in\cK \\
& \alpha_{ki} \ge 0 & \forall i\in[M], k\in\cK \\
& \overline{\beta}_{k\ell}, \underline{\beta}_{k\ell} \ge 0 & \forall k,\ell\in\cK \\
& \pi_{ei} \ge 0 & \forall e\in\cE, i\in[M] \\
& \overline{\rho}_{e\ell}, \underline{\rho}_{e\ell} \ge 0  & \forall e\in\cE, \ell\in\cK \\
& \xi_{kpi} \ge 0 & \forall k\in\cK,p\in \cP_k,i\in[M] \\
&\overline{\zeta}_{kp\ell}, \underline{\zeta}_{kp\ell}\ge 0 & \forall k,\ell\in\cK,p\in \cP_k.
\end{align*}

\section{Experiments}
\label{sec:experiments}

\subsection{Setup}

The aim of the experiments is to analyze the performance of solutions to the robust network design problem using discrete and polyhedral uncertainty sets,
respectively. We set all existing capacities $u_e$ to be zero, so that the effect of model choice becomes more visible. 

Note that an obvious additional comparison method is the kernel learning method from \cite{shang2017data} to construct polyhedral datasets. This method was implemented and tested. While it may have the advantage of giving a good representation of the data, this comes at the price of a large number of constraints and variables that define the uncertainty polyhedron. A single instance of the resulting linear program for the robust optimization problem was run for over 65 days without finding an optimal solution; hence, this method is not included in the comparison.

Using the data described in Section~\ref{sec:data}, we focus on the four months from beginning of May until end of August which provide the most complete
sets of data measurements. We base the training set on an arbitrarily chosen month, 07, which consists of 8,928 demand scenarios, but we removed outlier
scenarios as explained in Section~\ref{sec:data}, leaving us with $ T = 8,750 $ scenarios in the training set.

We calculate solutions based on the training set derived from month 07 measurements and then evaluate them on all the scenarios from the training set and
from the three other months (05, 06, 08), minimizing unsatisfied demand. Only the first-stage $\pmb{x}$-part of a solution is used for evaluation.
For discrete uncertainty, we calculate solutions based on clusterings with $K=100$ up to $K=8,600$ in steps of 100, and in addition using all $ T = 8,750 $
training scenarios (a total of 87 optimization problems and solutions). Clusters are calculated using the \texttt{kmeans} function of SciPy 1.2.1 under
Python 3.7.
For polyhedral uncertainty, we place hyperplanes using the \texttt{dual\_annealing} function from SciPy. We generate 140 hyperplanes this way. They are
collected in 28 polyhedra, where polyhedron $i$ uses all hyperplanes of polyhedron $i-1$, and five more in the order that they were generated. In total,
this means that 28 optimization problems with polyhedral uncertainty are solved.\footnote{All linear programs were solved using Cplex 12.8 on a virtual
Ubuntu server with ten Xeon CPU E7-2850 processors at 2.00 GHz speed and 23.5 GB RAM using only one core each.}

For a better comparison, the $87+28$ solutions are then also scaled up and down uniformly by multiplying the corresponding $\pmb{x}$ vector with a factor
$\lambda=0.5$ up to $\lambda=1.5$ (with step size 1/40). We consider these scaled versions as by construction, solutions based on
polyhedra are more conservative than those based on clusterings. By scaling solutions up and down, a more comprehensive comparison becomes possible.
Each of these $(87+28)\cdot 41$ solutions is then evaluated by calculating an optimal flow for each of the 8750 training scenarios and each of the $8928+8640+8928$ evaluation scenarios. In total, this means that over 166 million linear programs are solved for the evaluation. As there may not be sufficient capacity available to route all demand, we minimize the unsatisfied demand in each optimization problem. The corresponding model to evaluate solutions $\pmb{x}$ for a fixed scenario $\pmb{d}$ is hence as follows:
\begin{align*}
\min\ & \sum_{k\in \cK} h_k \\
\text{s.t. } & h_k \ge d_k - \sum_{p\in \cP_k} f_{kp} & \forall k\in\cK \\
& \sum_{k\in\cK} \sum_{p\in \cP_k:e\in p} f_{kp} \le u_e + x_e & \forall e\in\cE \\
& x_e \ge 0 & \forall e\in\cE \\
& f_{kp} \ge 0 & \forall k\in\cK, p\in \cP_k,
\end{align*}
where $h_k$ denotes the unsatisfied demand in commodity $k$. Note that the cost of a solution only depends on the choice of $\pmb{x}$. Additionally,
for each month, we calculate four measures to characterize the performance of a solution with regard to unsatisfied demand: average, {\cvs} (i.e., the average unsatisfied demand over the 25\% largest values), {\cvn}, and the maximum. These measures are chosen to reflect a broad spectrum of possible risk preferences of the decision maker.

\subsection{Results}

\begin{figure}[htbp]
\begin{center}
\begin{subfigure}{.4\textwidth}
\includegraphics[width=\linewidth]{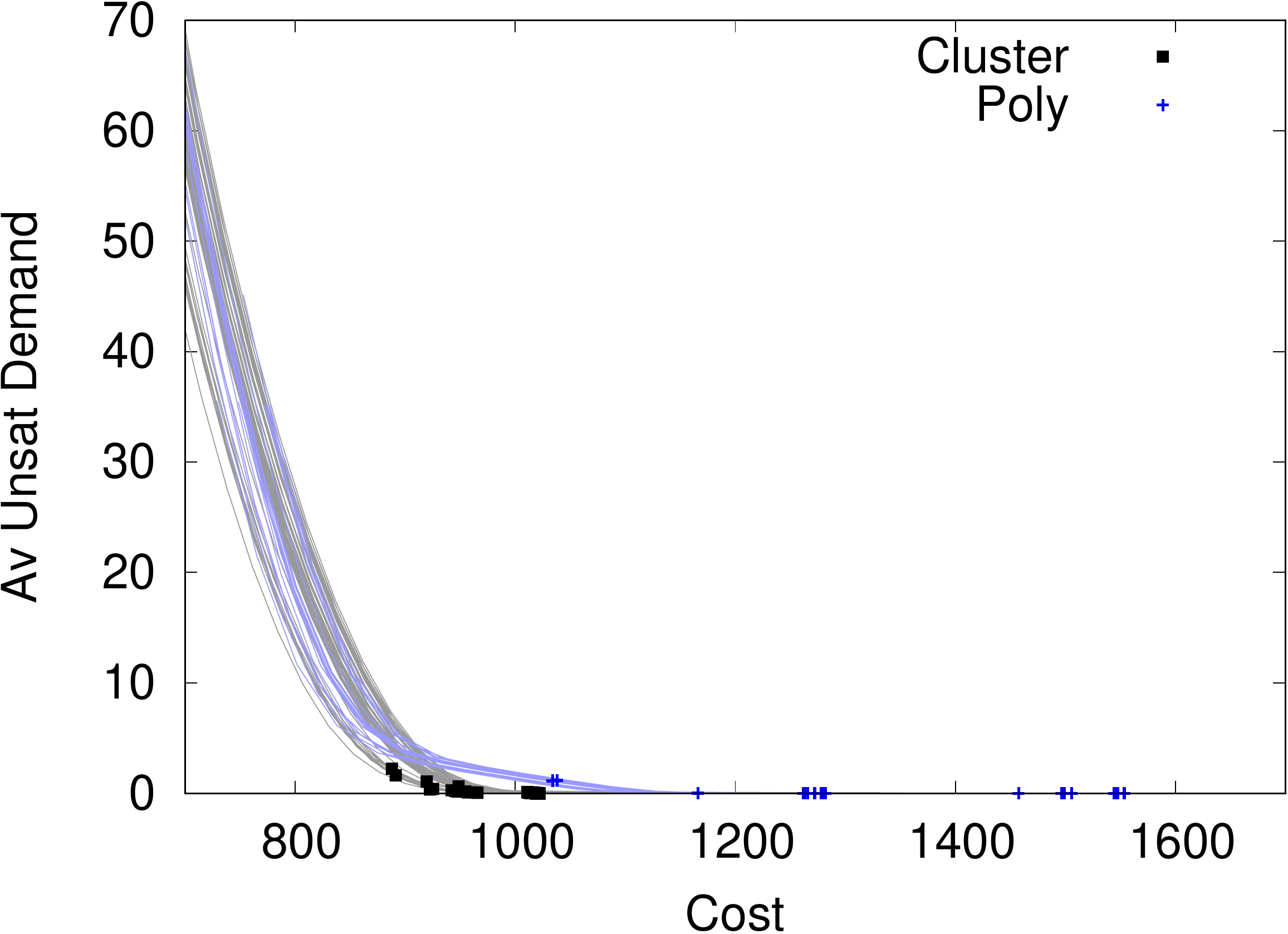}
\caption{Training set, average values.}\label{fig:data7av}
\end{subfigure}
\hspace{1cm}
\begin{subfigure}{.4\textwidth}
\includegraphics[width=\linewidth]{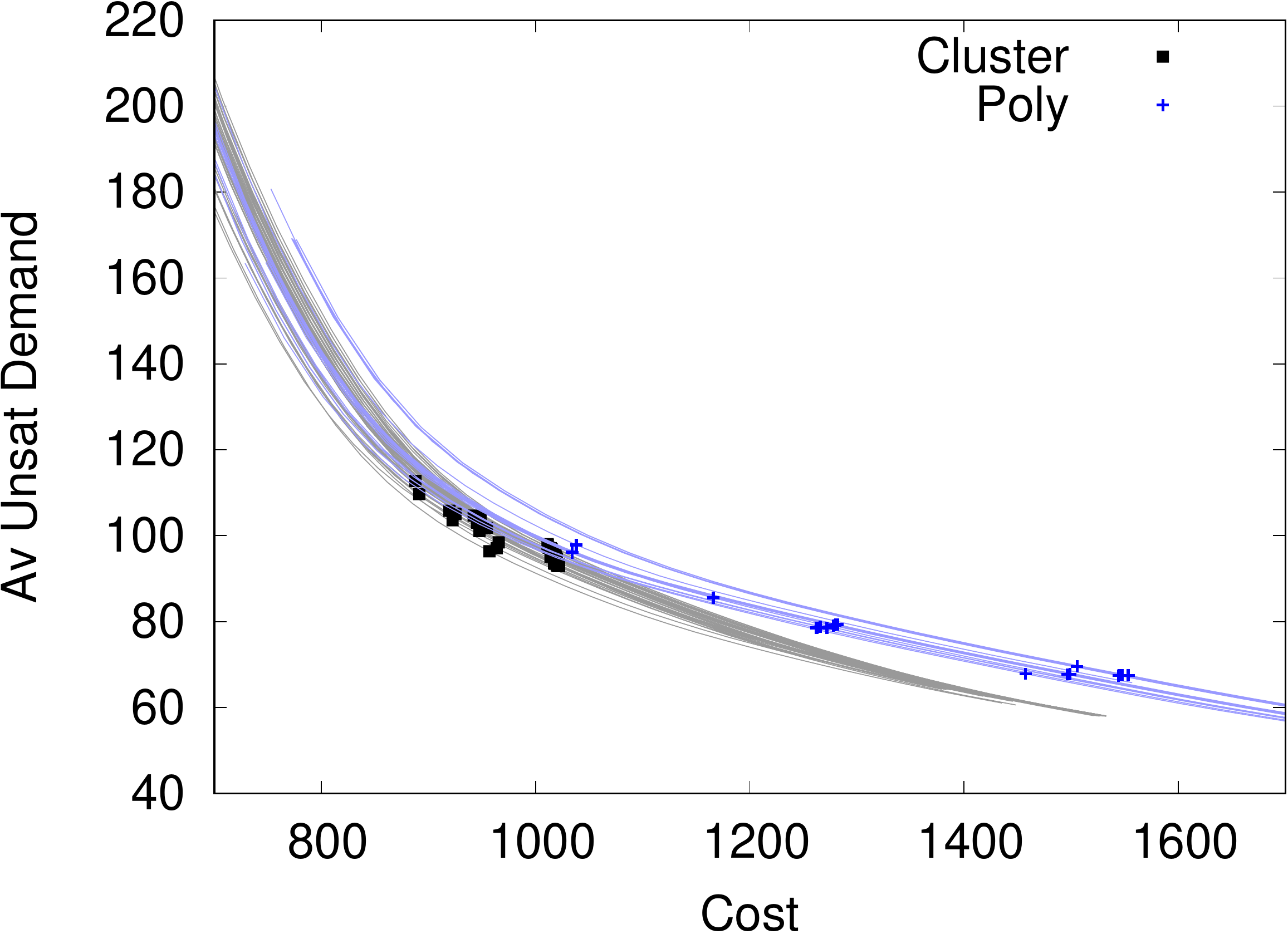}
\caption{Month 08, average values.}\label{fig:data8av}
\end{subfigure}

\vspace*{5mm}

\begin{subfigure}{.4\textwidth}
\includegraphics[width=\linewidth]{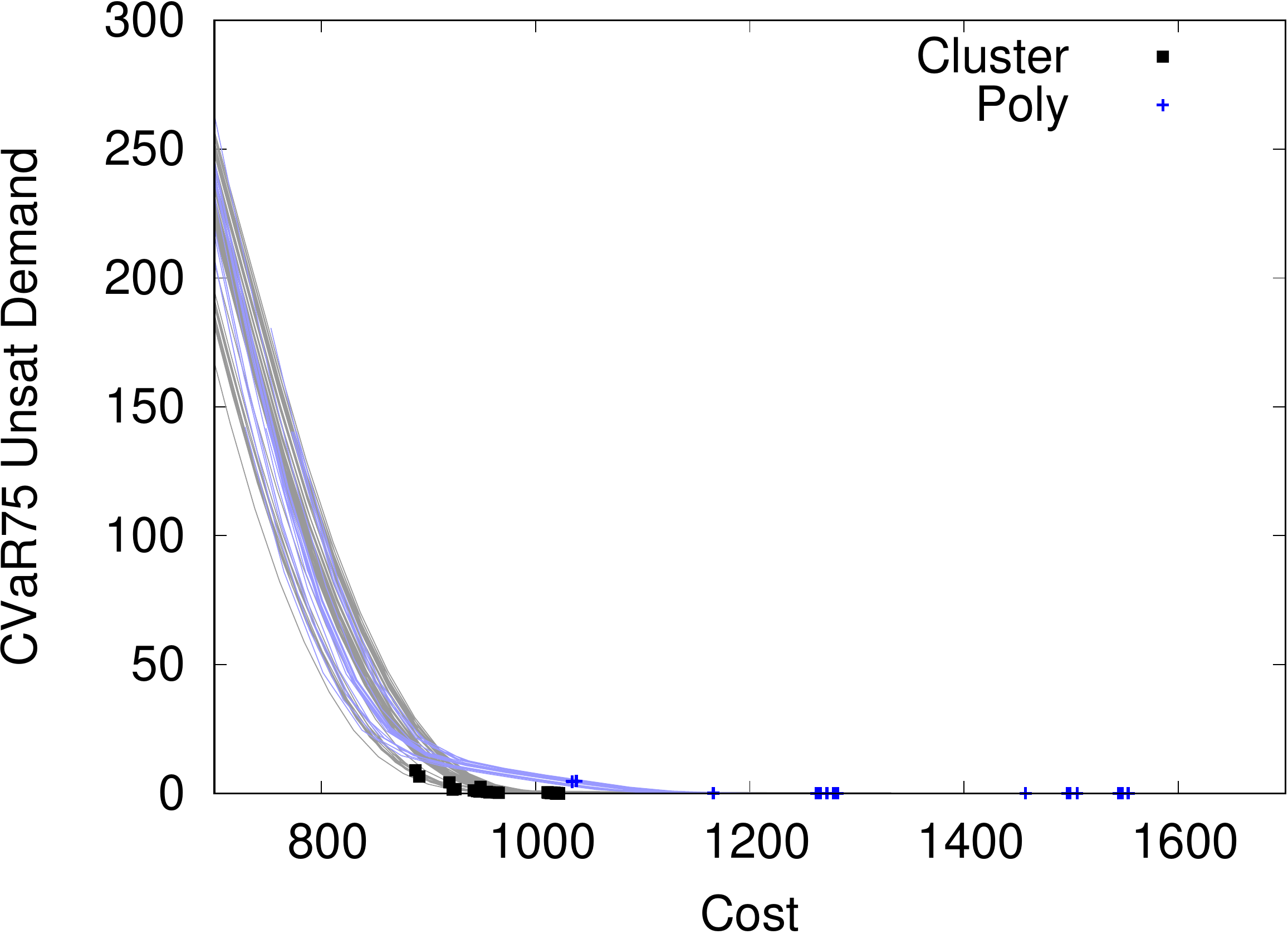}
\caption{Training set, {\cvs} values.}\label{fig:data7c75}
\end{subfigure}
\hspace{1cm}
\begin{subfigure}{.4\textwidth}
\includegraphics[width=\linewidth]{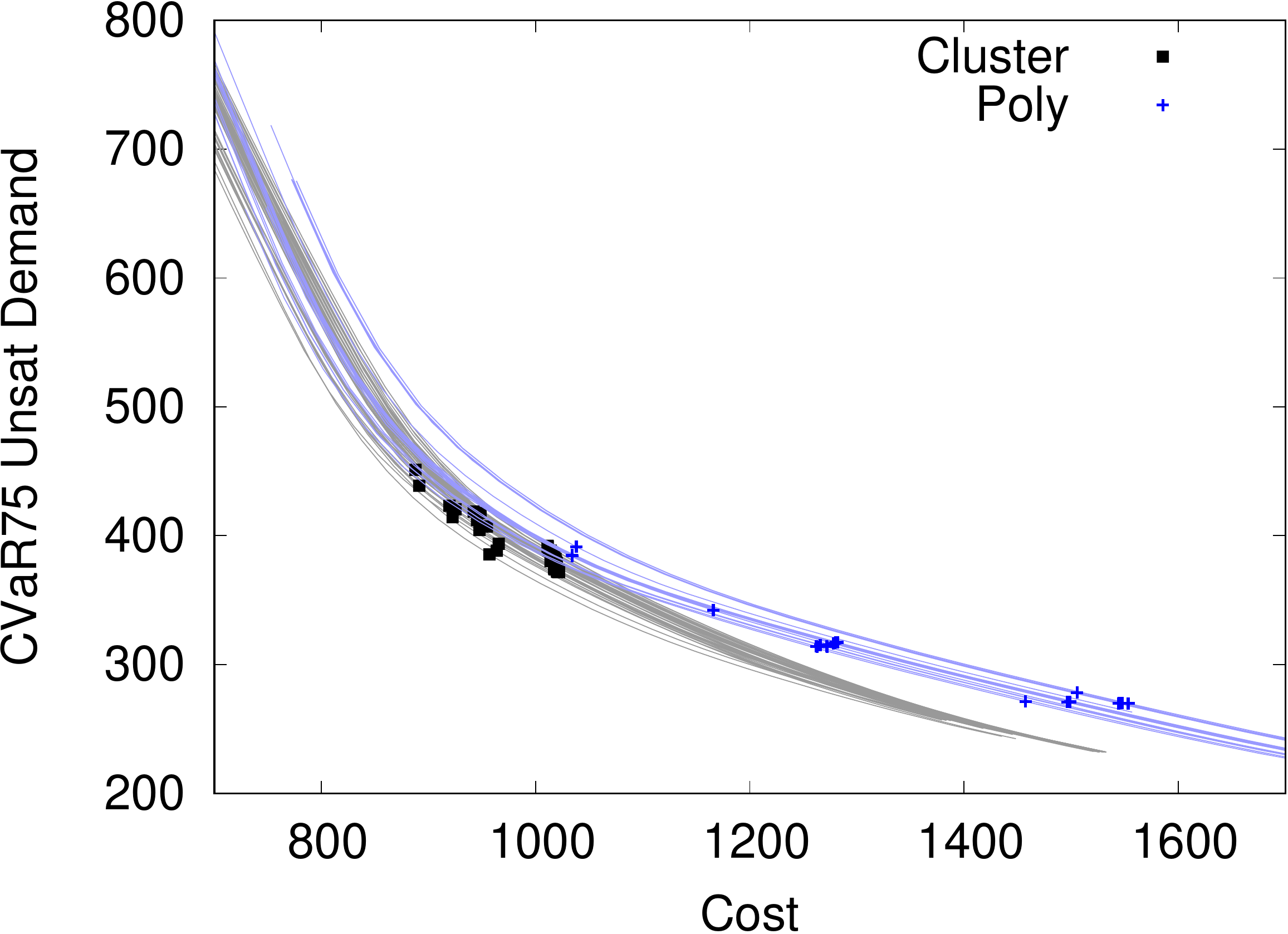}
\caption{Month 08, {\cvs} values.}\label{fig:data8c75}
\end{subfigure}

\vspace*{5mm}

\begin{subfigure}{.4\textwidth}
\includegraphics[width=\linewidth]{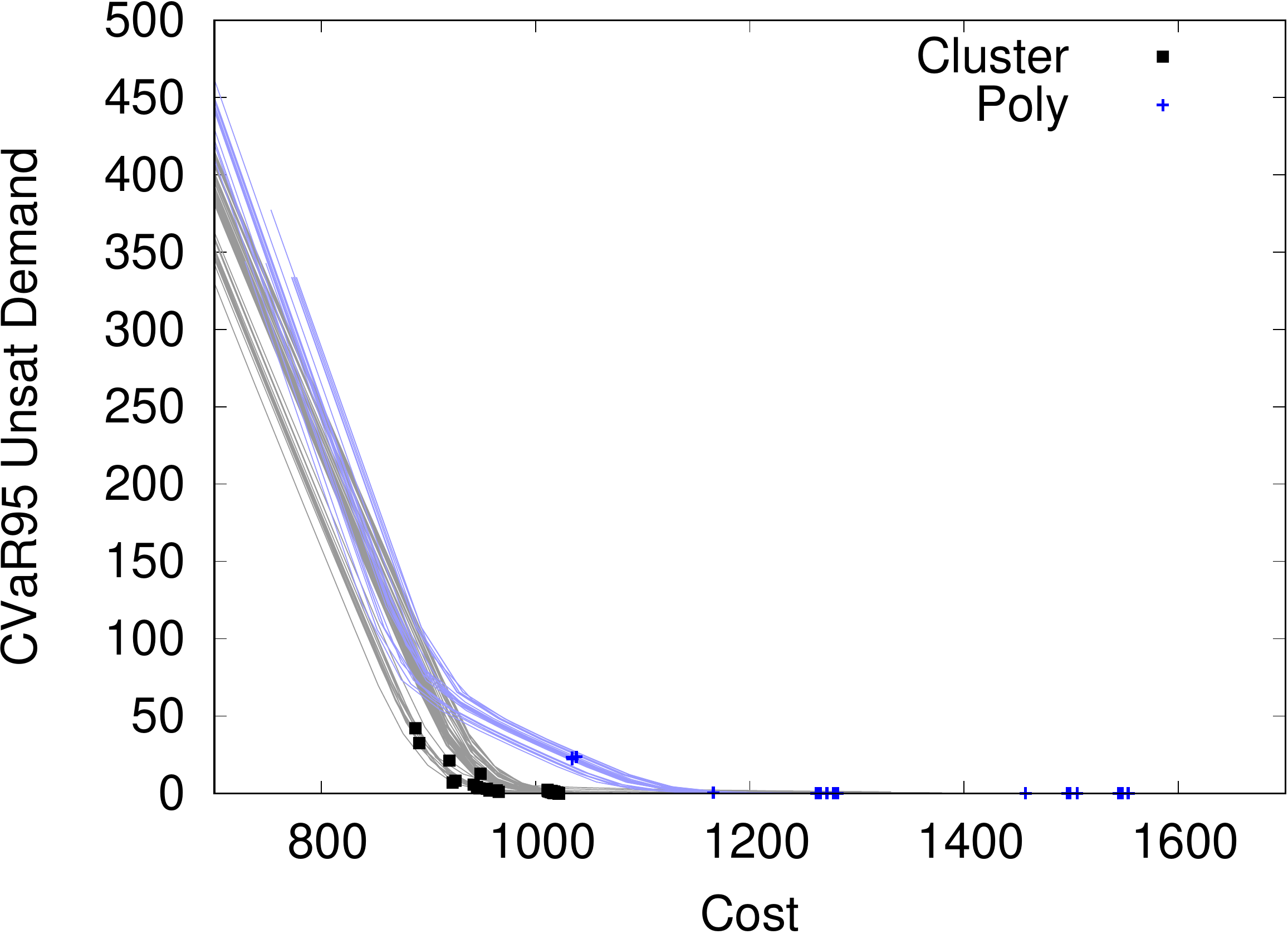}
\caption{Training set, {\cvn} values.}\label{fig:data7c95}
\end{subfigure}
\hspace{1cm}
\begin{subfigure}{.4\textwidth}
\includegraphics[width=\linewidth]{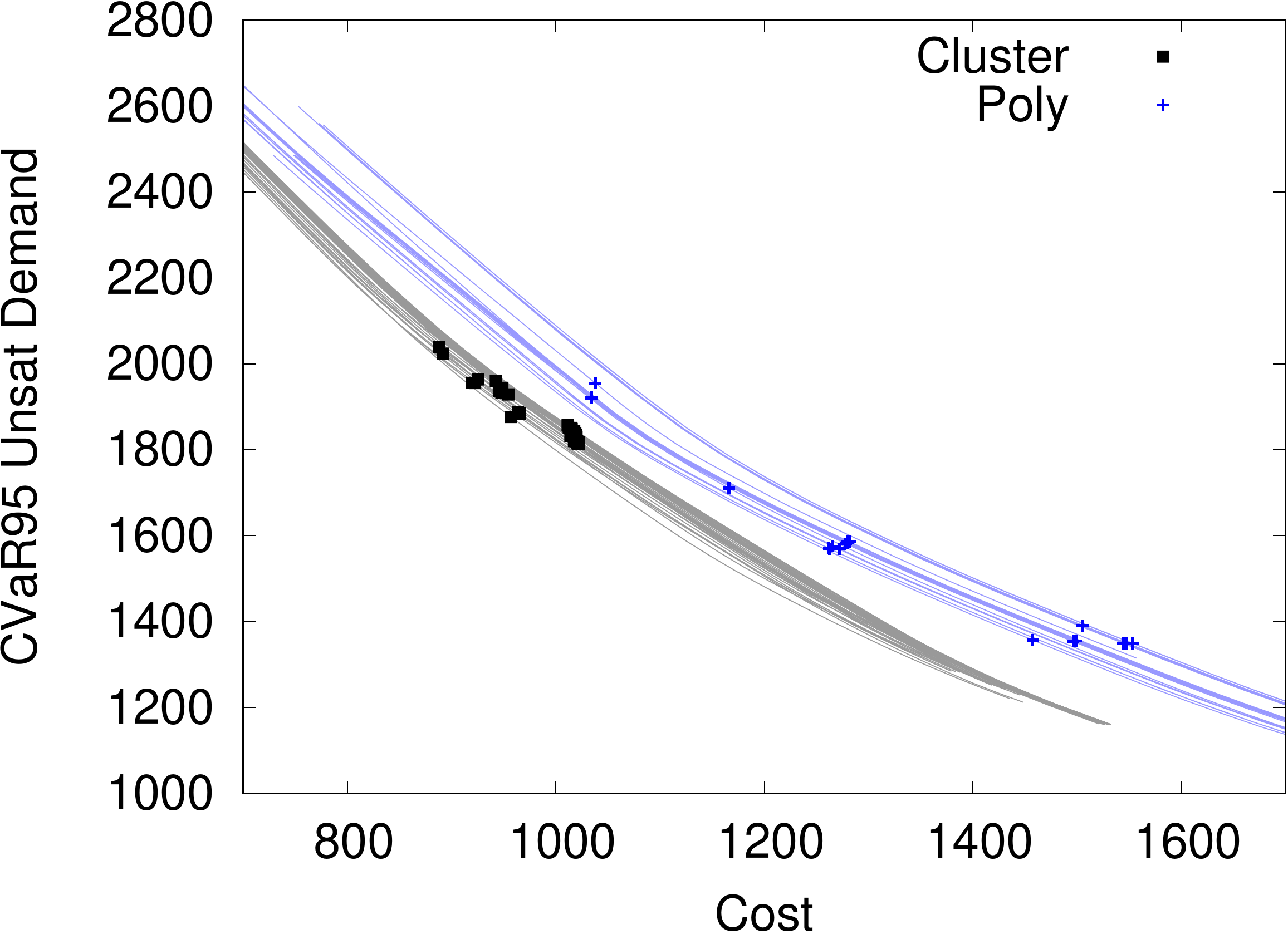}
\caption{Month 08, {\cvn} values.}\label{fig:data8c95}
\end{subfigure}

\vspace*{5mm}

\begin{subfigure}{.4\textwidth}
\includegraphics[width=\linewidth]{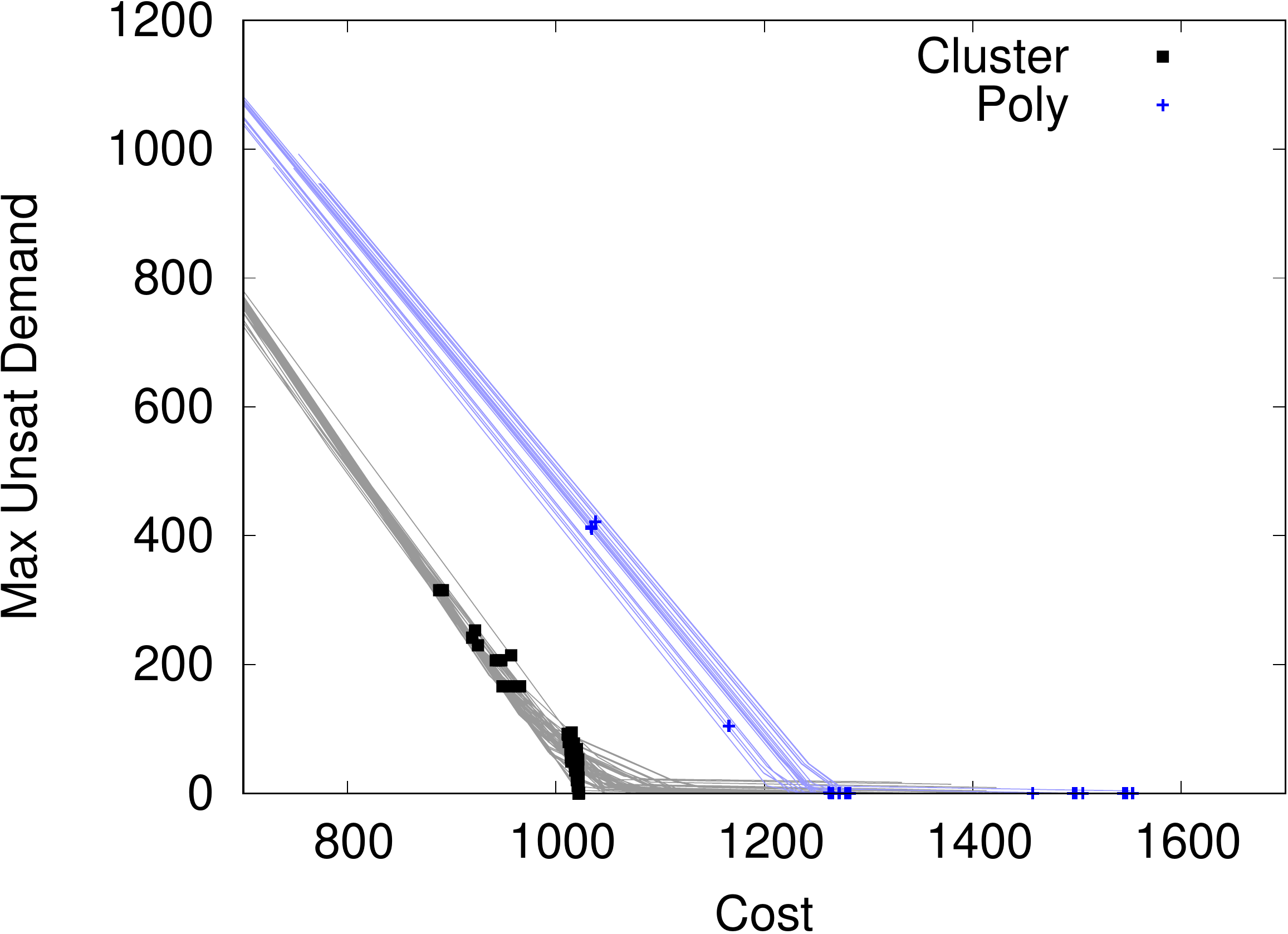}
\caption{Training set, maximum values.}\label{fig:data7max}
\end{subfigure}
\hspace{1cm}
\begin{subfigure}{.4\textwidth}
\includegraphics[width=\linewidth]{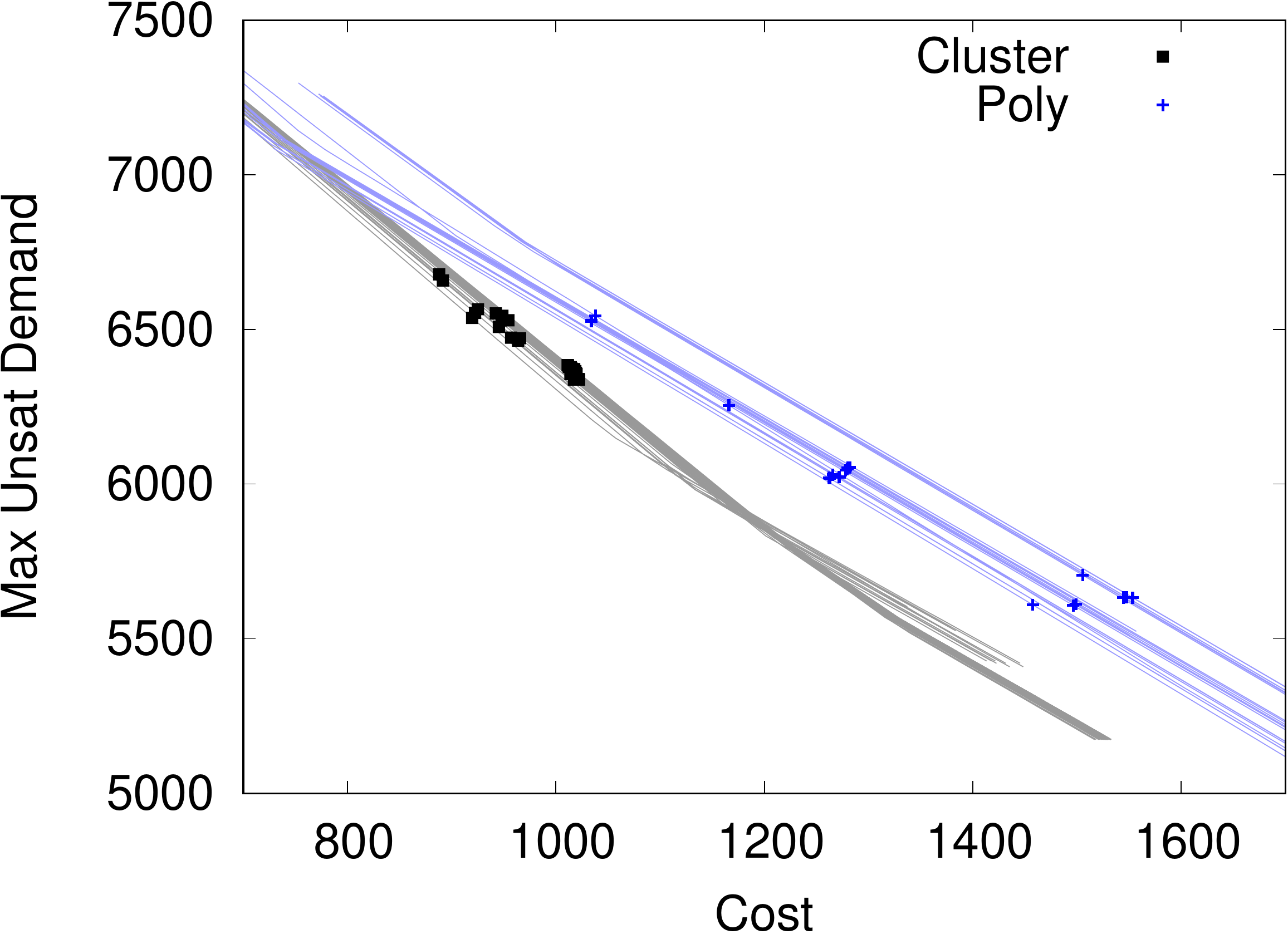}
\caption{Month 08, maximum values.}\label{fig:data8max}
\end{subfigure}
\caption{Results for training set and month 08.}\label{res78}
\end{center}
\end{figure}

We first discuss the performance of solutions on the training set (month 07) in the left column of Figure~\ref{res78}. On the horizontal axis, we show the
costs of solutions, while the vertical axis shows the four measures of unsatisfied demand. Each point corresponds to a solution (87 black squares
corresponding to the discrete uncertainty solutions, 28 blue crosses corresponding to the polyhedral uncertainty solutions). The lines show the
performance of the scaled solutions.

Consider Figure~\ref{fig:data7max}. By construction, we know that the discrete uncertainty solution with $\cU^d_{8750}$ has zero unsatisfied demand on the
training set, and is the cheapest possible solution to do so. Most polyhedral solutions use conservative outer approximations of the training data and thus
also have zero unmet demand, but at higher costs. We can also see that solutions that use fewer clusters or more hyperplanes become less
conservative, allowing unsatisfied demand at lower solution costs. The behavior we see in Figure~\ref{fig:data7max} is to be expected by design. The open
question is whether it can also be observed on evaluation data.

Figures~\ref{fig:data7av}, \ref{fig:data7c75} and \ref{fig:data7c95} show the average, {\cvs}, and {\cvn} performance on the training set, respectively. Here the differences between both types of solution are much less pronounced; we see that blue and black lines overlap, indicating a similar performance of solution types.

Compare this performance to the right-hand column of Figure~\ref{res78}, where the performance on month 08 is presented. The order of magnitude of
unsatisfied demand has increased for each type of solution: whereas in Figure~\ref{fig:data7max} we can reach zero unsatisfied demand, the same solutions
have between five and seven thousand units of maximum unmet demand in Figure~\ref{fig:data8max}. But the relative performance between the solution types is
similar. Whereas solutions based on polyhedral uncertainty generally have a higher degree of robustness at higher investment costs, it is possible to scale
solutions based on clustered data up to reach solutions with a similar degree of robustness at lower costs. This is particularly visible for the high
risk-adverse measures in Figures~\ref{fig:data8c95} and \ref{fig:data8max}, whereas these performance differences are less clear-cut for the less
risk-adverse measures in Figures~\ref{fig:data8av} and \ref{fig:data8c75}.

Figure~\ref{res56} in Appendix~\ref{sec:app} shows the results for months 05 and 06, where the same observations apply as for month 08.

In terms of solution quality, i.e., trade-off between investment costs and unsatisfied demand, we thus find the following result: Solutions based on
discrete uncertainty outperform solutions based on polyhedral uncertainty when using high risk-adverse
performance metrics such as maximum or {\cvn} of unsatisfied demand. This is less clear for less risk-adverse metrics such as expected value or {\cvs} of
unsatisfied demand, but the clustering approach is still superior in most cases.

\begin{figure}[htb]
\begin{center}
\includegraphics[width=0.5\textwidth]{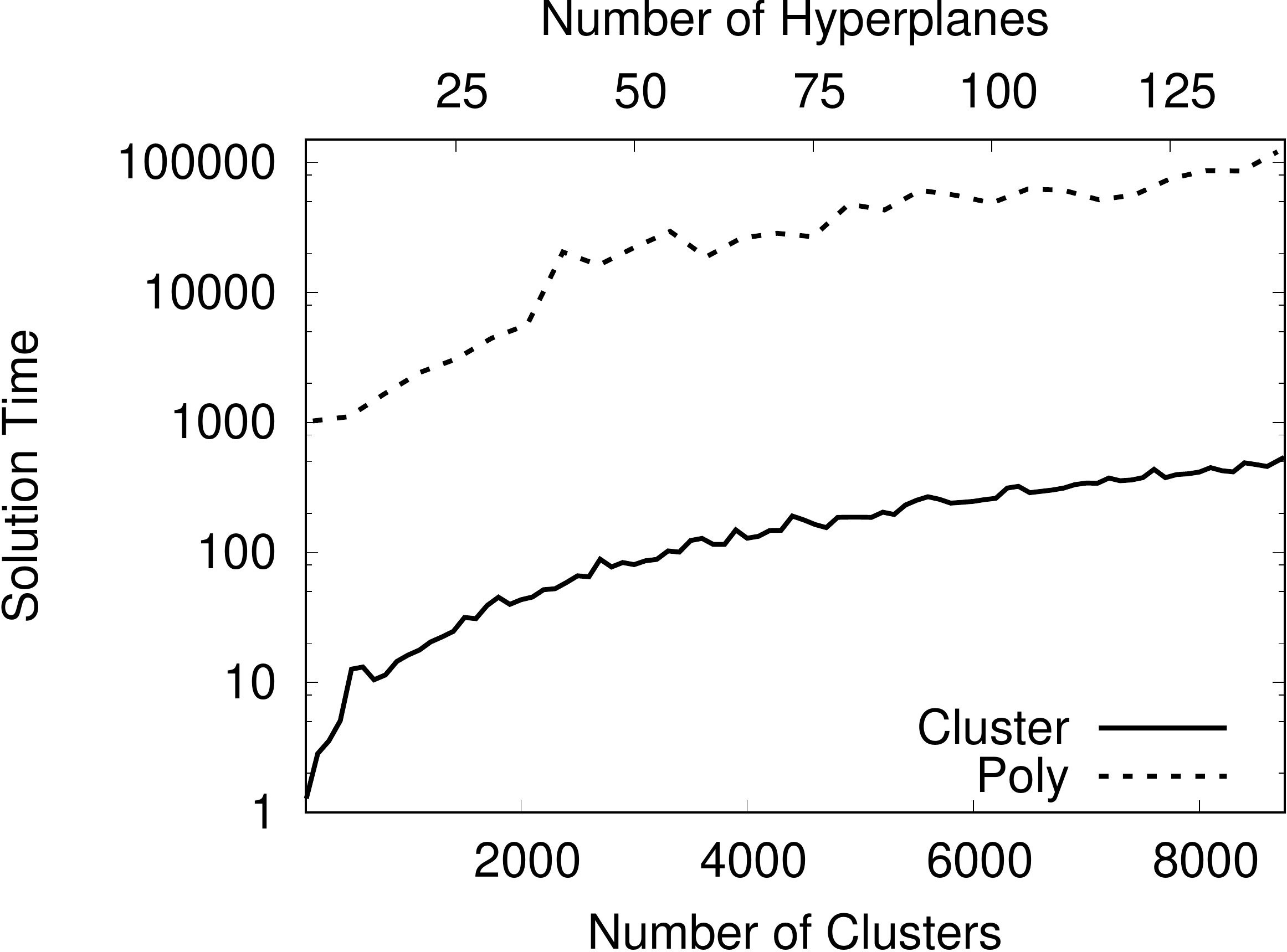}
\caption{Solution times for discrete and polyhedral uncertainty.}\label{fig:time}
\end{center}
\end{figure}

We now consider the time required to solve the corresponding robust optimization models. Figure~\ref{fig:time} shows the Cplex solution time for discrete
and for polyhedral uncertainty, which depends on the size of the uncertainty set (note the two different horizontal axes and the logarithmic vertical
axis).

It can be seen that even the largest discrete model (that uses all training scenarios directly) is still easier to solve than the smallest polyhedral model (using five hyperplanes in addition to the lower and upper bounds). So this experiment reveals that using discrete uncertainty sets not only results in a better solution quality, they are also easier to solve.

\section{Conclusions}
\label{sec:conclusions}

In the robust optimization literature, frequently both discrete and polyhedral uncertainty sets are being used. In this paper we compared the resulting
solutions using real-world data for a network expansion problem. We describe how to construct uncertainty sets based on clustering the training data and by separating training data from noise by placing hyperplanes. In our computational
study we found that solutions based on discrete uncertainty models outperform solutions based on polyhedral models in most performance metrics and are
also easier to compute. The strong performance of discrete uncertainty sets is in line with evidence from the experiment on shortest path data performed in
\cite{Chassein2018}. This also indicates that the current network design literature, which has a strong focus on polyhedral models, may benefit from
considering simple discrete models more.


%

\appendix

\section{Additional Experimental Results}
\label{sec:app}

\begin{figure}[htbp]
\begin{center}
\begin{subfigure}{.4\textwidth}
\includegraphics[width=\linewidth]{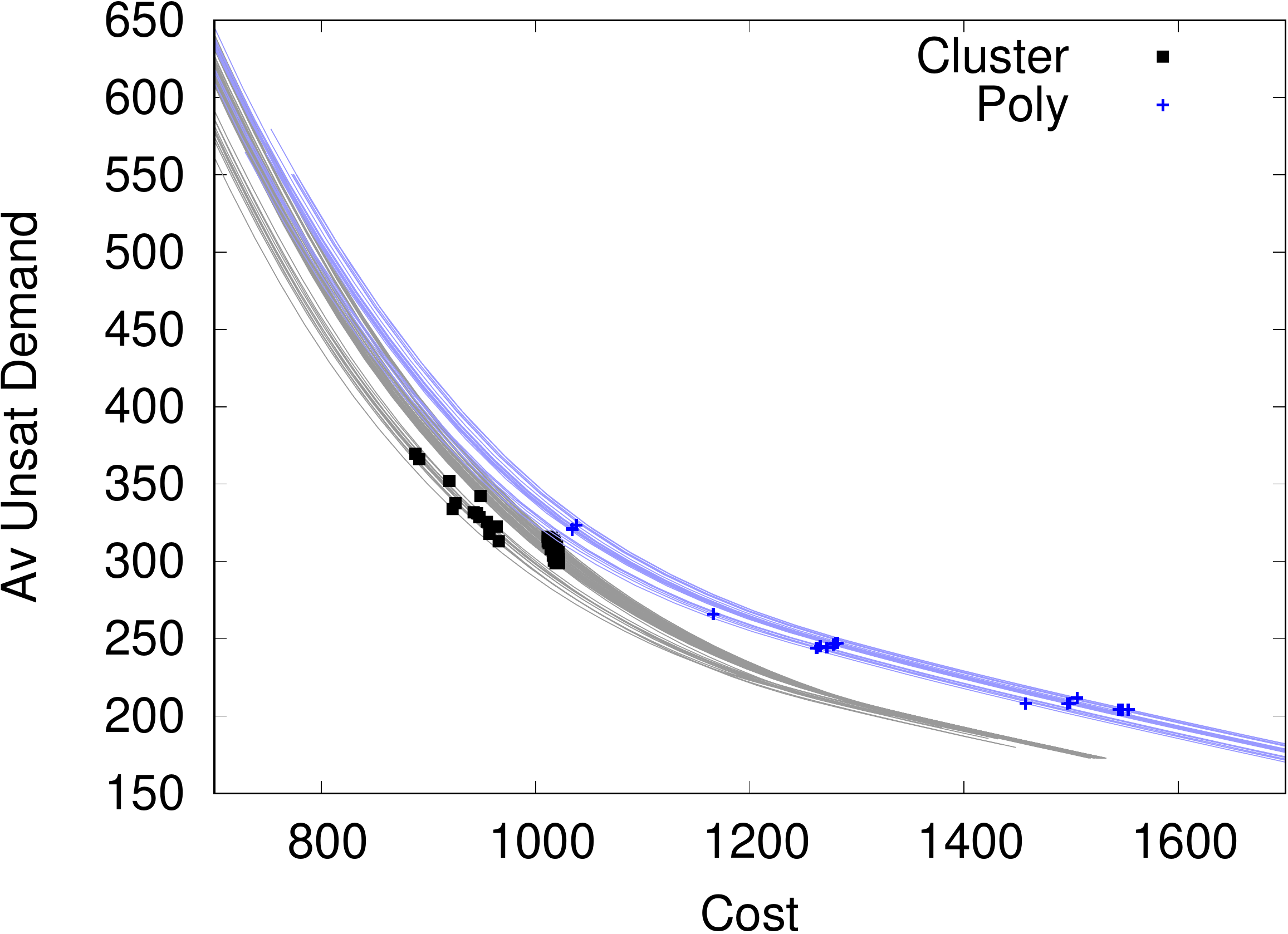}
\caption{Month 05, average values.}\label{fig:data5av}
\end{subfigure}
\hspace{1cm}
\begin{subfigure}{.4\textwidth}
\includegraphics[width=\linewidth]{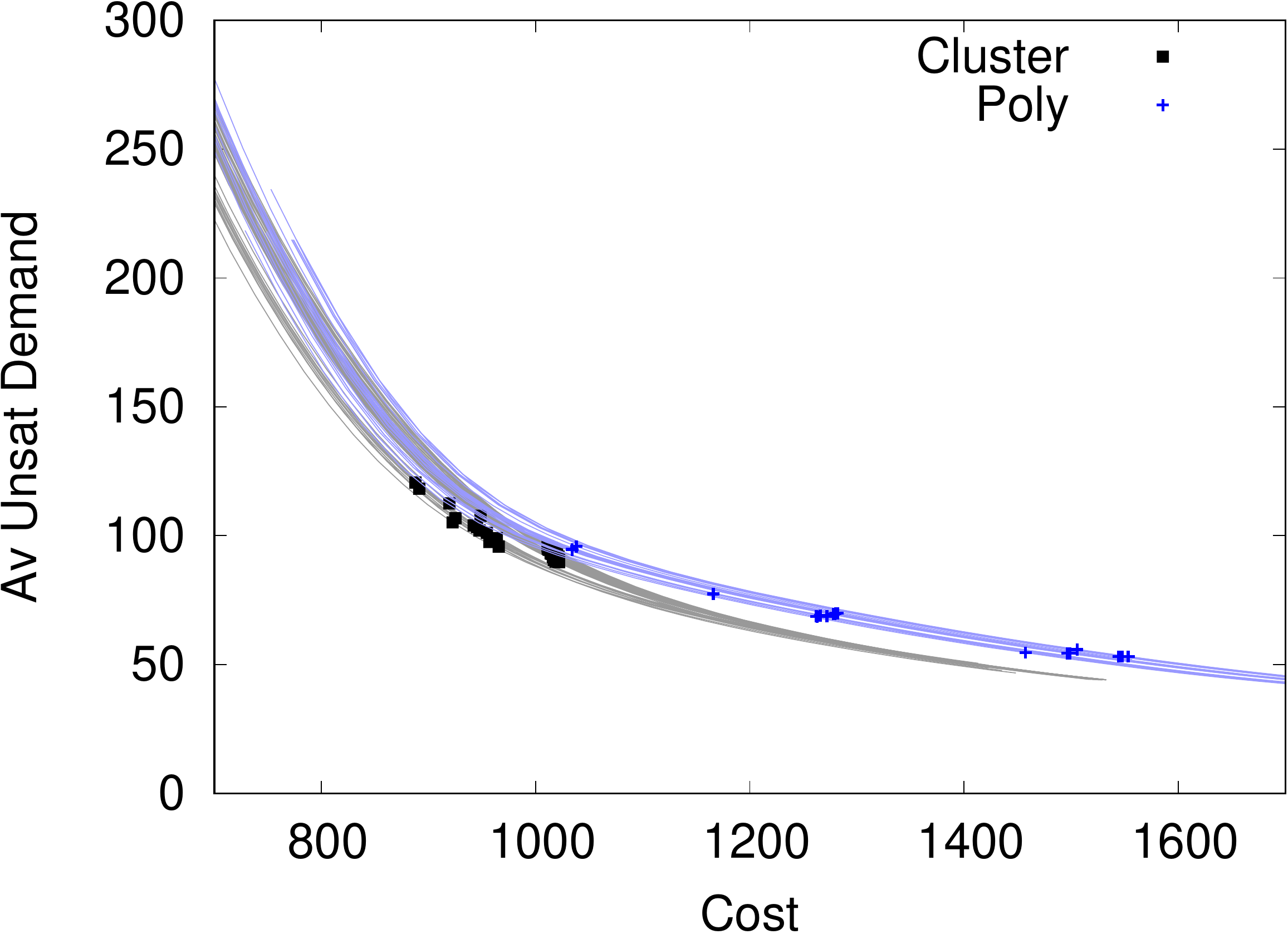}
\caption{Month 06, average values.}\label{fig:data6av}
\end{subfigure}

\vspace*{5mm}

\begin{subfigure}{.4\textwidth}
\includegraphics[width=\linewidth]{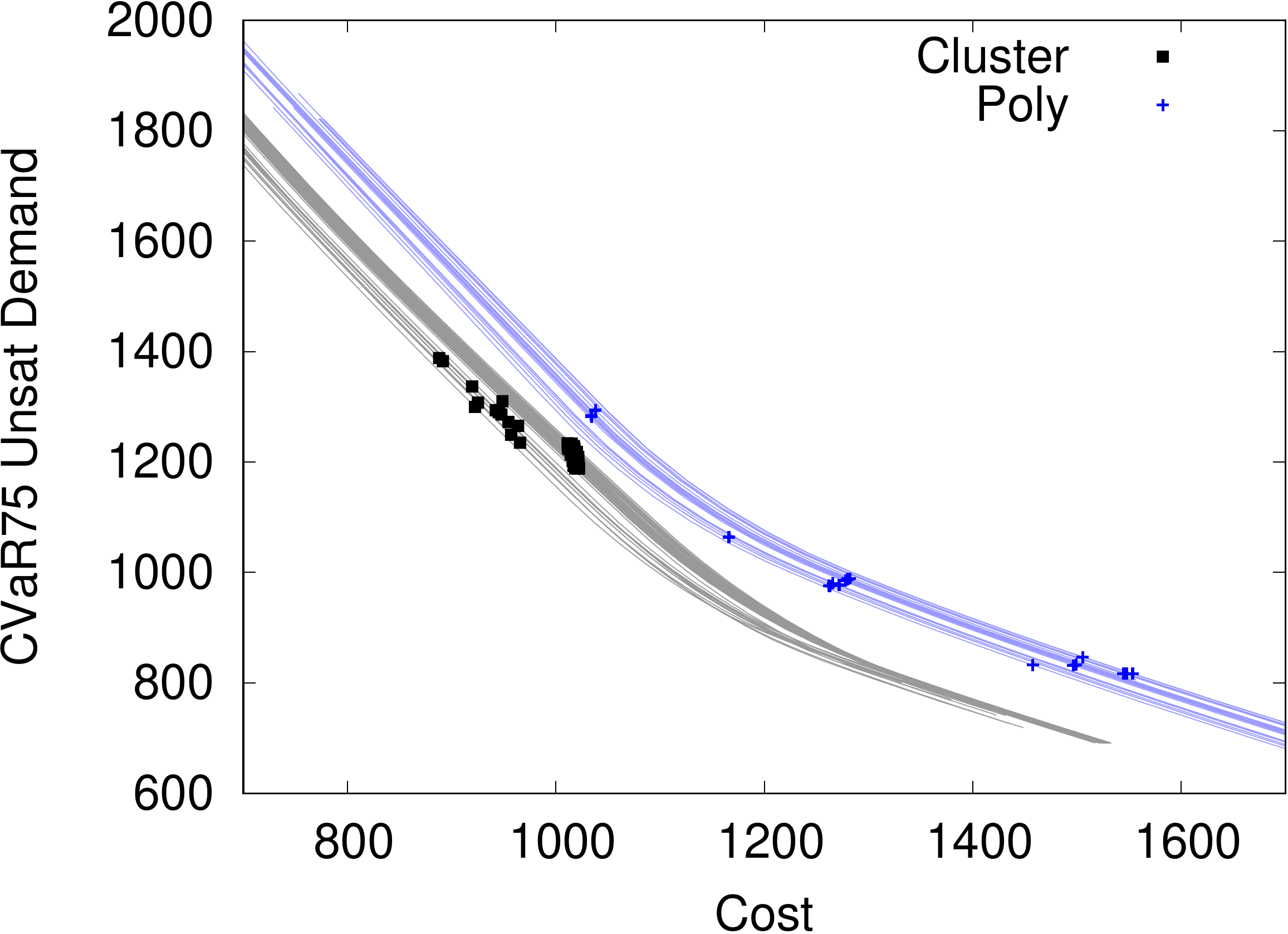}
\caption{Month 05, {\cvs} values.}\label{fig:data5c75}
\end{subfigure}
\hspace{1cm}
\begin{subfigure}{.4\textwidth}
\includegraphics[width=\linewidth]{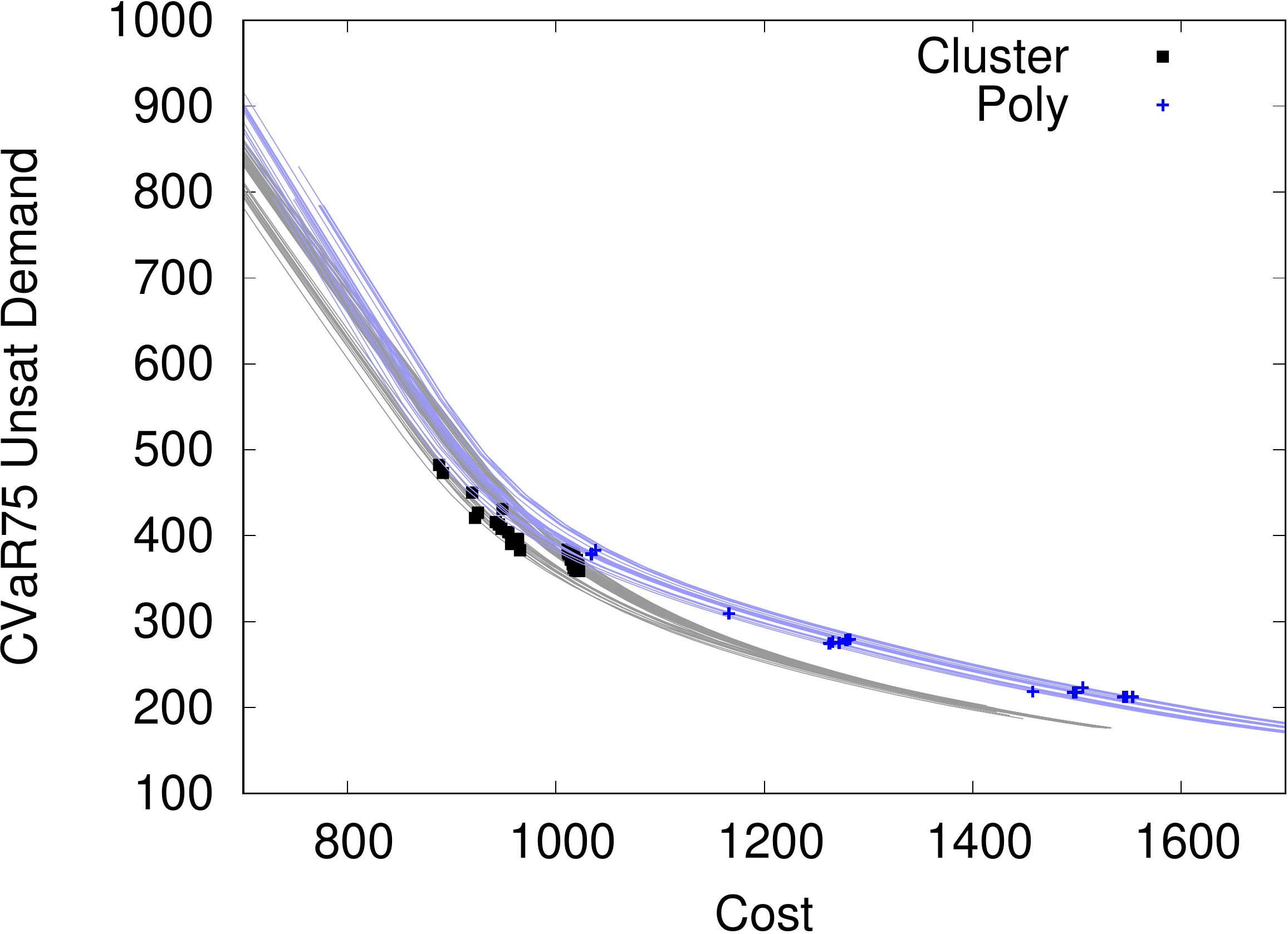}
\caption{Month 06, {\cvs} values.}\label{fig:data6c75}
\end{subfigure}

\vspace*{5mm}

\begin{subfigure}{.4\textwidth}
\includegraphics[width=\linewidth]{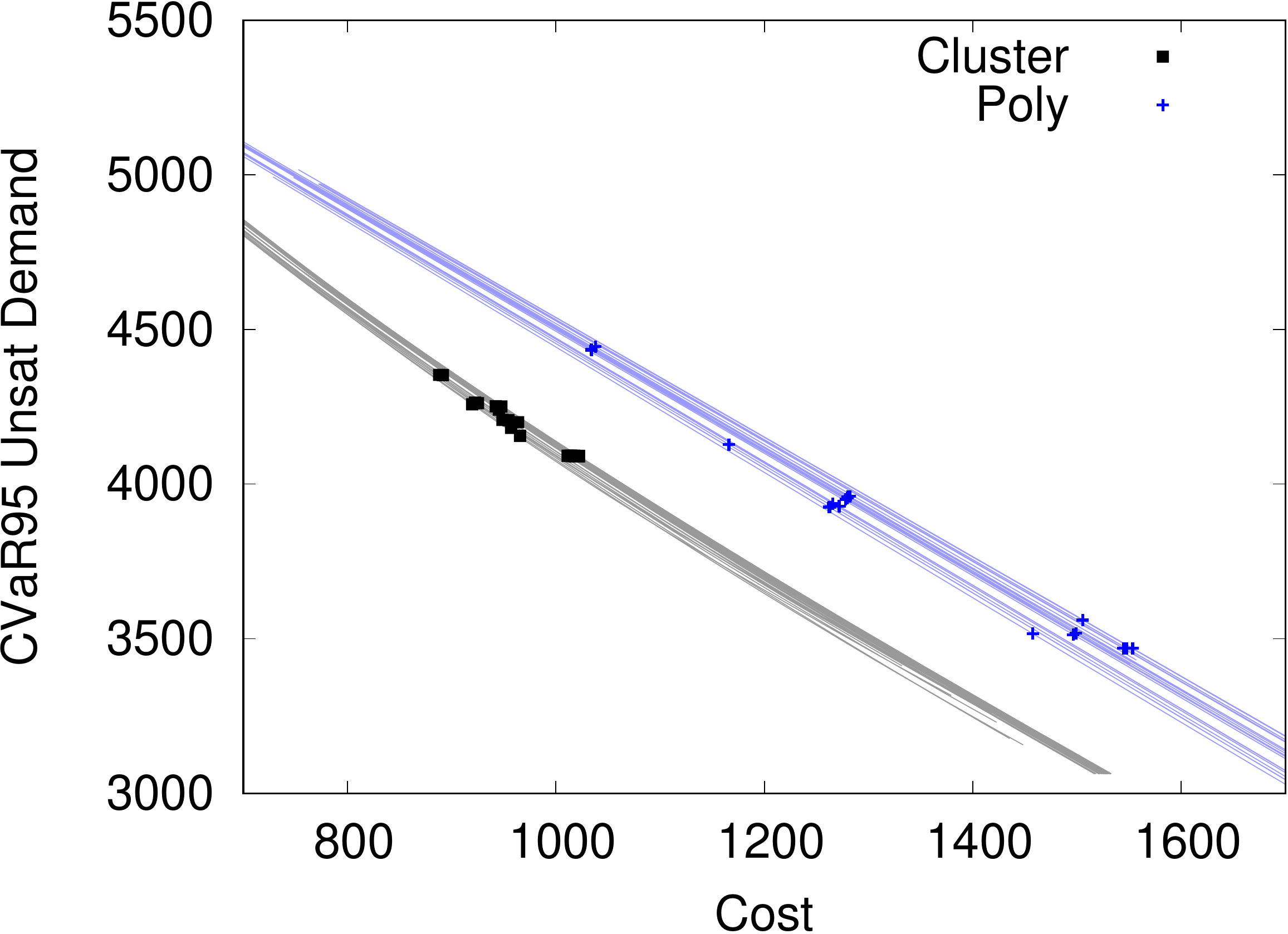}
\caption{Month 05, {\cvn} values.}\label{fig:data5c95}
\end{subfigure}
\hspace{1cm}
\begin{subfigure}{.4\textwidth}
\includegraphics[width=\linewidth]{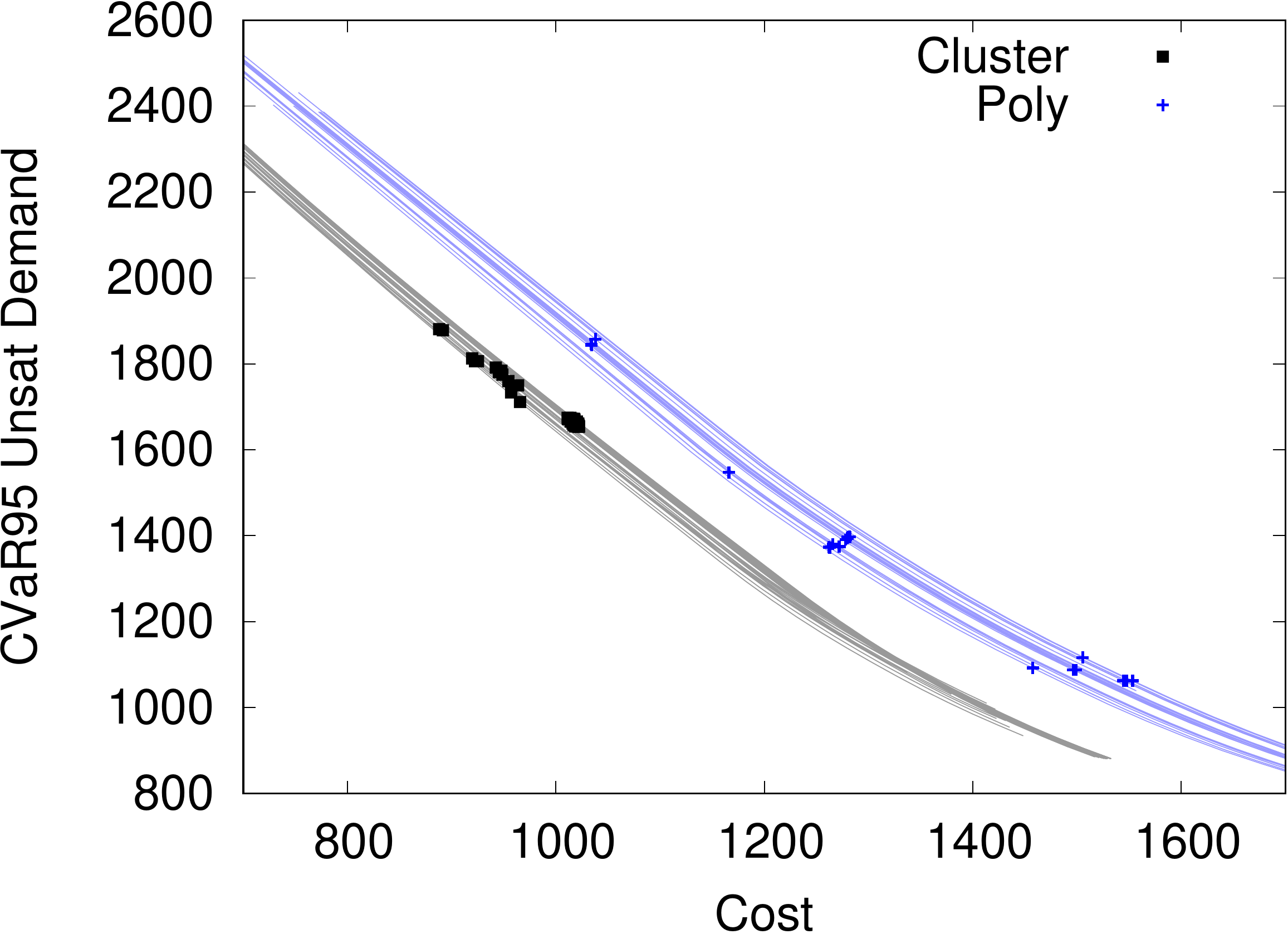}
\caption{Month 06, {\cvn} values.}\label{fig:data6c95}
\end{subfigure}

\vspace*{5mm}

\begin{subfigure}{.4\textwidth}
\includegraphics[width=\linewidth]{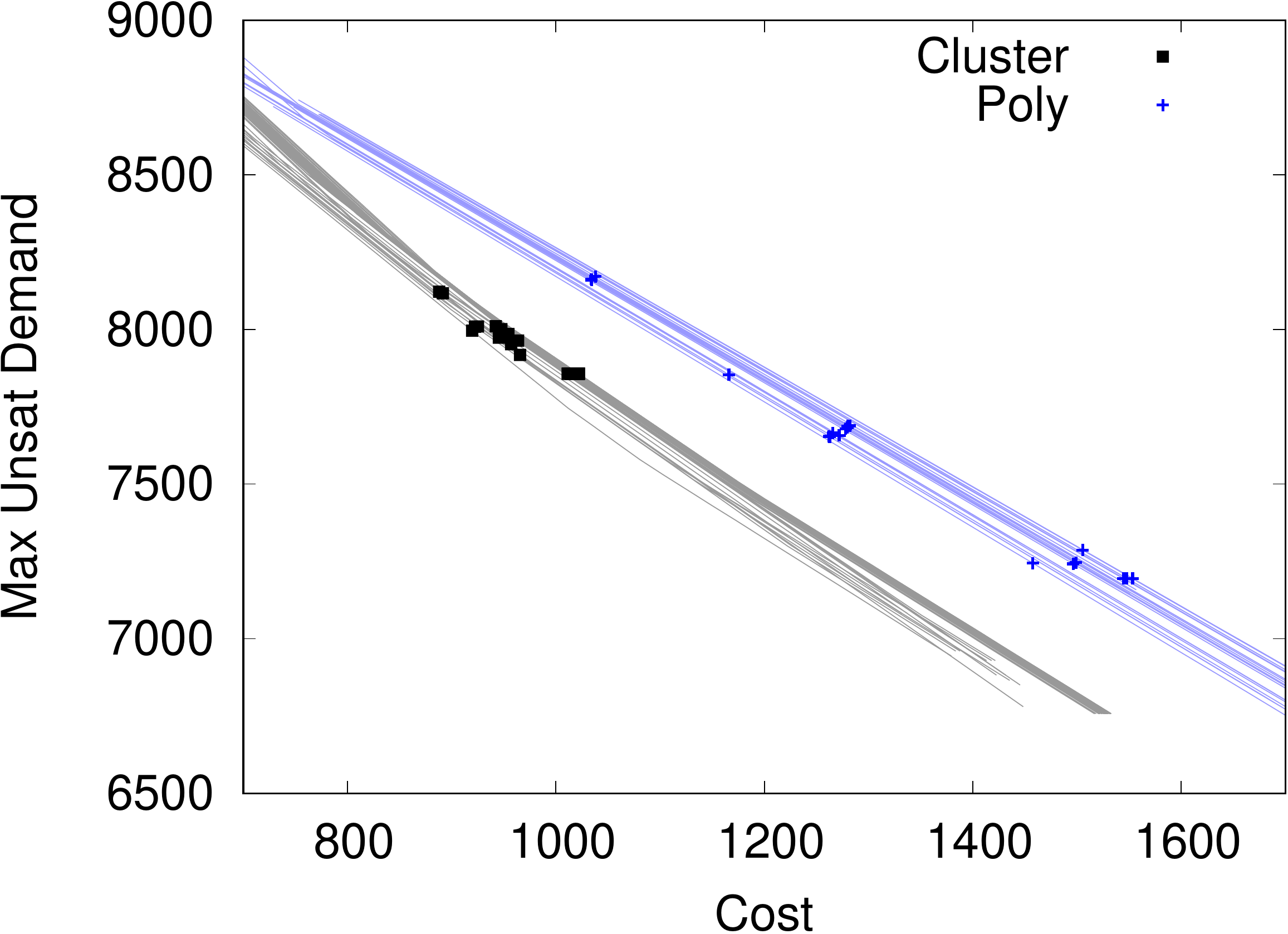}
\caption{Month 05, maximum values.}\label{fig:data5max}
\end{subfigure}
\hspace{1cm}
\begin{subfigure}{.4\textwidth}
\includegraphics[width=\linewidth]{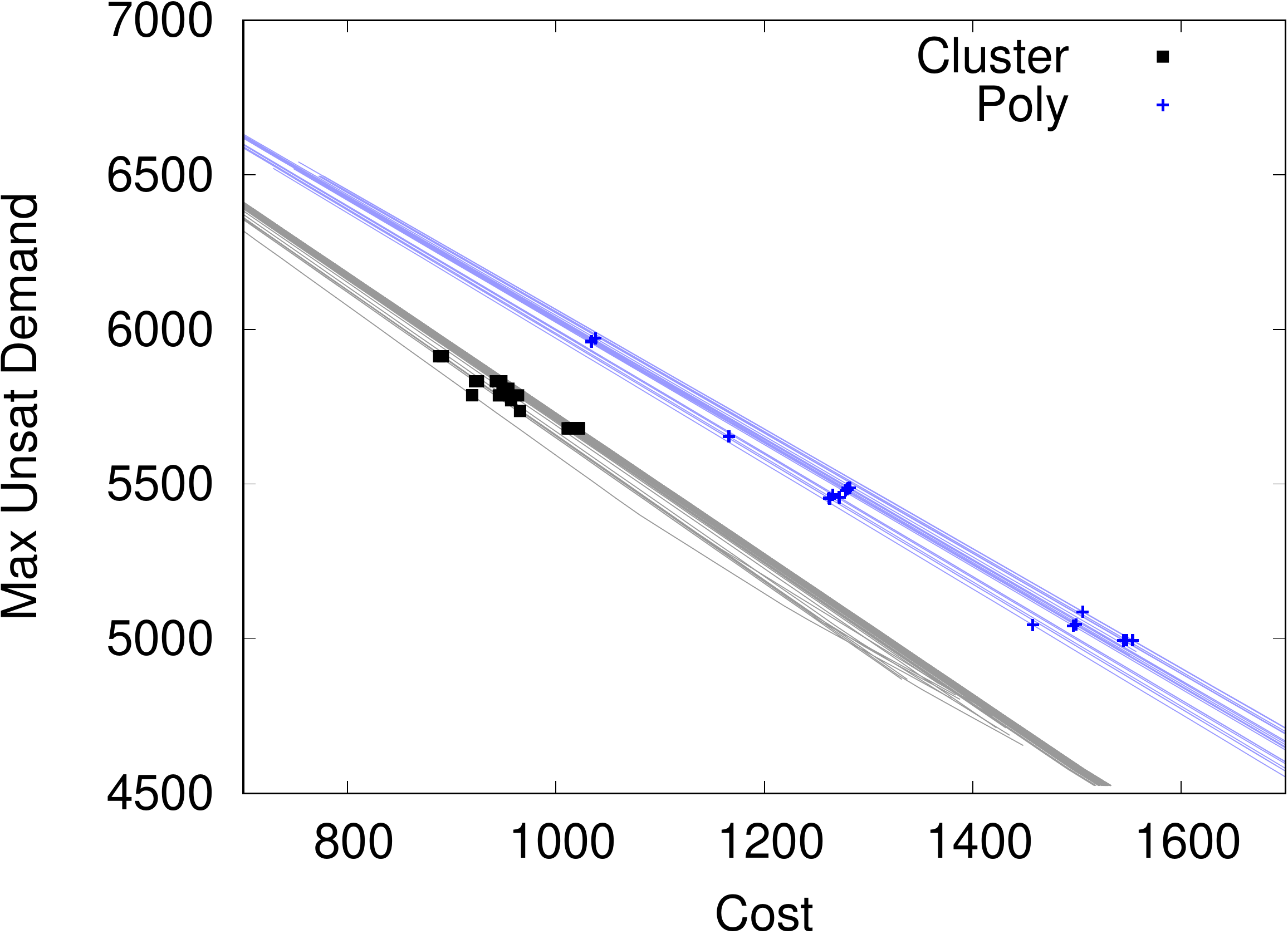}
\caption{Month 06, maximum values.}\label{fig:data6max}
\end{subfigure}
\caption{Results for months 05 and 06.}\label{res56}
\end{center}
\end{figure}


\begin{thebibliography}{10}

\bibitem{Atamt_rk_2007}
A.~Atamt{\"u}rk and M.~Zhang.
\newblock Two-stage robust network flow and design under demand uncertainty.
\newblock {\em Operations Research}, 55(4):662--673, 2007.

\bibitem{Babonneau2013}
F.~Babonneau, J.-P. Vial, O.~Klopfenstein, and A.~Ouorou.
\newblock Robust capacity assignment solutions for telecommunications networks
  with uncertain demands.
\newblock {\em Networks}, 62(4):255--272, 2013.

\bibitem{Ben-Tal2009}
A.~Ben-Tal, L.~El~Ghaoui, and A.~Nemirovski.
\newblock {\em Robust Optimization}.
\newblock Princeton University Press, 2009.

\bibitem{Ben_Tal_2004}
A.~Ben-Tal, A.~Goryashko, E.~Guslitzer, and A.~Nemirovski.
\newblock Adjustable robust solutions of uncertain linear programs.
\newblock {\em Mathematical programming}, 99(2):351--376, 2004.

\bibitem{Bertsekas1998}
D.~P. Bertsekas.
\newblock {\em Network Optimization: Continuous and Discrete Models}.
\newblock Athena Scientific, Belmont, Massachusetts, 1998.

\bibitem{Bertsimas2011}
D.~Bertsimas, D.~B. Brown, and C.~Caramanis.
\newblock Theory and applications of robust optimization.
\newblock {\em SIAM review}, 53(3):464--501, 2011.

\bibitem{Bertsimas2017}
D.~Bertsimas, V.~Gupta, and N.~Kallus.
\newblock Data-driven robust optimization.
\newblock {\em Mathematical Programming}, 167(2):235--292, 2018.

\bibitem{Bertsimas_2004}
D.~Bertsimas and M.~Sim.
\newblock The price of robustness.
\newblock {\em Operations research}, 52(1):35--53, 2004.

\bibitem{campbell2015bayesian}
T.~Campbell and J.~P. How.
\newblock Bayesian nonparametric set construction for robust optimization.
\newblock In {\em 2015 American Control Conference (ACC)}, pages 4216--4221.
  IEEE, 2015.

\bibitem{Chassein2018}
A.~Chassein, T.~Dokka, and M.~Goerigk.
\newblock Algorithms and uncertainty sets for data-driven robust shortest path
  problems.
\newblock {\em European Journal of Operational Research}, 274(2):671 -- 686,
  2019.

\bibitem{chassein2018approximating}
A.~Chassein, M.~Goerigk, A.~Kasperski, and P.~Zieli{\'n}ski.
\newblock Approximating combinatorial optimization problems with the ordered
  weighted averaging criterion.
\newblock {\em European Journal of Operational Research}, 286(3):828--838,
  2020.

\bibitem{Duffield1999}
N.~G. Duffield, P.~Goyal, A.~Greenberg, P.~Mishra, K.~K. Ramakrishnan, and
  J.~E. van~der Merive.
\newblock A flexible model for resource management in virtual private networks.
\newblock In {\em Proceedings of the conference on Applications, technologies,
  architectures, and protocols for computer communication}, pages 95--108,
  1999.

\bibitem{Fingerhut_1997}
J.~A. Fingerhut, S.~Suri, and J.~S. Turner.
\newblock Designing least-cost nonblocking broadband networks.
\newblock {\em Journal of Algorithms}, 24(2):287--309, 1997.

\bibitem{Gendron1999}
B.~Gendron, T.~G. Crainic, and A.~Frangioni.
\newblock Multicommodity capacitated network design.
\newblock In {\em Telecommunications network planning}, pages 1--19. Springer,
  1999.

\bibitem{goerigk2016algorithm}
M.~Goerigk and A.~Sch{\"o}bel.
\newblock Algorithm engineering in robust optimization.
\newblock In {\em Algorithm engineering}, pages 245--279. Springer
  International Publishing, 2016.

\bibitem{Koster_2013}
A.~M. C.~A. Koster, M.~Kutschka, and C.~Raack.
\newblock Robust network design: Formulations, valid inequalities, and
  computations.
\newblock {\em Networks}, 61(2):128--149, 2013.

\bibitem{Magnanti1984}
T.~L. Magnanti and R.~T. Wong.
\newblock Network design and transportation planning: Models and algorithms.
\newblock {\em Transportation science}, 18(1):1--55, 1984.

\bibitem{Minoux1989}
M.~Minoux.
\newblock Networks synthesis and optimum network design problems: Models,
  solution methods and applications.
\newblock {\em Networks}, 19(3):313--360, 1989.

\bibitem{ning2017data}
C.~Ning and F.~You.
\newblock Data-driven adaptive nested robust optimization: general modeling
  framework and efficient computational algorithm for decision making under
  uncertainty.
\newblock {\em AIChE Journal}, 63(9):3790--3817, 2017.

\bibitem{ning2018data}
C.~Ning and F.~You.
\newblock Data-driven decision making under uncertainty integrating robust
  optimization with principal component analysis and kernel smoothing methods.
\newblock {\em Computers \& Chemical Engineering}, 112:190--210, 2018.

\bibitem{Ord_ez_2007}
F.~Ord{\'{o}}{\~{n}}ez and J.~Zhao.
\newblock Robust capacity expansion of network flows.
\newblock {\em Networks}, 50(2):136--145, 2007.

\bibitem{Orlowski2010}
S.~Orlowski, R.~Wess\"{a}ly, M.~Pi\'{o}ro, and A.~Tomaszewski.
\newblock {SNDlib} 1.0-{S}urvivable {N}etwork {D}esign {L}ibrary.
\newblock {\em Networks}, 55(3):276--286, 2010.

\bibitem{Ouorou2007}
A.~Ouorou and J.-P. Vial.
\newblock A model for robust capacity planning for telecommunications networks
  under demand uncertainty.
\newblock In {\em 6th International Workshop on Design and Reliable
  Communication Networks (DRCN 2007)}, pages 1--4. IEEE, 2007.

\bibitem{Pessoa_2015}
A.~A. Pessoa and M.~Poss.
\newblock Robust network design with uncertain outsourcing cost.
\newblock {\em INFORMS Journal on Computing}, 27(3):507--524, 2015.

\bibitem{Poss_2012}
M.~Poss and C.~Raack.
\newblock Affine recourse for the robust network design problem: Between static
  and dynamic routing.
\newblock {\em Networks}, 61(2):180--198, 2013.

\bibitem{shang2017data}
C.~Shang, X.~Huang, and F.~You.
\newblock Data-driven robust optimization based on kernel learning.
\newblock {\em Computers \& Chemical Engineering}, 106:464--479, 2017.

\end{thebibliography}
\end{document}